\documentclass[a4paper,11pt,reqno]{amsart}

\usepackage{amssymb}
\usepackage{color}

\newtheorem{theorem}{Theorem}
\newtheorem{lemma}[theorem]{Lemma}

\newtheorem{proposition}[theorem]{Proposition}
\numberwithin{theorem}{section}
\numberwithin{equation}{section}

\theoremstyle{definition}
\newtheorem{definition}[theorem]{Definition}

\newtheorem{remark}[theorem]{Remark}

\newtheorem*{remark*}{Remark}

\newcommand{\Z}{{\mathbb Z}}

\newcommand{\cA}{{\mathcal A}}

\newcommand{\cC}{{\mathcal C}}

\newcommand{\J}{{\mathcal J}}
\newcommand{\cJ}{{\mathcal J}}
\newcommand{\cH}{{\mathcal H}}

\newcommand{\cL}{{\mathcal L}}
\newcommand{\cM}{{\mathcal M}}

\newcommand{\cS}{{\mathcal S}}
\newcommand{\cT}{{\mathcal T}}
\newcommand{\cQ}{{\mathcal Q}}
\newcommand{\cO}{{\mathcal O}}

\newcommand{\cU}{{\mathcal U}}
\newcommand{\g}{{\mathfrak g}}
\newcommand{\gl}{{\mathfrak{gl}}}

\newcommand{\h}{{\mathfrak h}}
\newcommand{\m}{{\mathfrak m}}

\newcommand{\frso}{{\mathfrak{so}}}  
\newcommand{\so}{{\mathfrak{so}}}
\newcommand{\frsp}{{\mathfrak{sp}}}

\newcommand{\spl}{{\mathfrak{sl}}}

\newcommand{\anti}{{\mathfrak{skew}}}
\newcommand{\sym}{{\mathfrak{sym}}}

 \DeclareMathOperator{\eespan}{span}
 \DeclareMathOperator{\End}{End}
 \DeclareMathOperator{\Hom}{Hom}

 \DeclareMathOperator{\tr}{tr}
 \DeclareMathOperator{\Ad}{Ad}
  
 \DeclareMathOperator{\ad}{ad}
 \DeclareMathOperator{\Der}{Der}
\DeclareMathOperator{\Inder}{Inder}
 \DeclareMathOperator{\Mat}{Mat}

 \def\dimens{\mathop{\rm dim}\nolimits}

\begin{document}

\title[Irreducible Lie-Yamaguti algebras of Generic Type]%
{Irreducible Lie-Yamaguti algebras of\\
Generic Type}

\author{Pilar Benito}

\thanks{Supported by the
Spanish Ministerio de Educaci\'on y Ciencia and FEDER (MTM
2007-67884-C04-02,03). Pilar Benito and Fabi\'an Mart\'{\i}n-Herce
also acknowledges support from the  Comunidad Aut\'onoma de La Rioja
(ANGI2005/05,06),
 and Alberto Elduque from
the Diputaci\'on General de Arag\'on (Grupo de Investigaci\'on de
\'Algebra).}

\address{Departamento de Matem\'aticas y Computaci\'on y Centro de Investigaci\'on de Inform\'atica, Matem\'aticas y Estad\'{\i}stica, Universidad de
La Rioja, 26004 Logro\~no, Spain}

\email{pilar.benito@unirioja.es}



\author{Alberto Elduque}

\address{Departamento de Matem\'aticas e
Instituto Universitario de Matem\'aticas y Aplicaciones, Universidad de
Zaragoza, 50009 Zaragoza, Spain}

\email{elduque@unizar.es}

\author{Fabi\'an Mart\'{\i}n-Herce}

\address{Departamento de Matem\'aticas y Computaci\'on, Universidad de
La Rioja, 26004 Logro\~no, Spain}

\email{fabian.martin@dmc.unirioja.es}

\date{July 20, 2009}



\begin{abstract}
Lie-Yamaguti algebras (or generalized Lie triple systems) are
binary-ternary algebras intimately related to reductive homogeneous
spaces. The Lie-Yamaguti algebras which are irreducible as modules
over their inner derivation algebras are the algebraic
counterpart of the isotropy irreducible homogeneous spaces.

These systems splits into three disjoint types: adjoint type,
non-simple type and generic type. The systems of the first two types
were classified in a previous paper through a generalized Tits
Construction of Lie algebras. In this paper, the Lie-Yamaguti
algebras of generic type are classified by relating them to several
other nonassociative algebraic systems: Lie and Jordan algebras and
triple systems, Jordan pairs or Freudenthal triple systems.
\end{abstract}

\maketitle


\section{Introduction}



Let $G$ be a connected Lie group with Lie algebra $\g$, $H$ a closed
subgroup of $G$, and let $\h$ be the associated subalgebra of $\g$.
The corresponding homogeneous space $M=G/H$ is said to be
\emph{reductive} (\cite[\S 7]{Nom54}) in case there is a subspace
$\m$ of $\g$ such that $\g=\h\oplus\m$ and $\Ad(H)(\m)\subseteq \m$.

In this situation, Nomizu proved \cite[Theorem 8.1]{Nom54} that
there is a one-to-one correspondence between the set of all
$G$-invariant affine connections on $M$ and the set of bilinear
multiplications $\alpha:\m\times\m\rightarrow\m$ such that the
restriction of $\Ad(H)$ to $\m$ is a subgroup of the automorphism
group of the nonassociative algebra $(\m,\alpha)$.

There exist natural binary and ternary products defined in $\m$,
given by
\begin{equation}\label{eq:binter}
\begin{split}
&x\cdot y = \pi_{\m}\bigl([x,y]\bigr),\\
&[x,y,z]=\bigl[\pi_{\h}([x,y]),z],
\end{split}
\end{equation}
for any $x,y,z\in\m$, where $\pi_{\h}$ and $\pi_{\m}$ denote the
projections on $\h$ and $\m$ respectively, relative to the reductive
decomposition $\g=\h\oplus\m$. Note that the condition
$\Ad(H)(\m)\subseteq \m$ implies the condition $[\h,\m]\subseteq
\m$, the converse being valid if $H$ is connected.

There are two distinguished invariant affine connections: the
 natural connection (or canonical connection of the first kind),
which corresponds to the bilinear multiplication given by
$\alpha(x,y)=\frac{1}{2} x\cdot y$ for any $x,y\in\m$, which has
trivial torsion, and the canonical connection corresponding to the
trivial multiplication: $\alpha(x,y)=0$ for any $x,y\in\m$. In case
the reductive homogeneous space is symmetric, so $[\m,\m]\subseteq
\h$, these two connections coincide. For the canonical connection,
the torsion and curvature tensors are given on the tangent space to
the point $eH\in M$ ($e$ denotes the identity element of $G$), which
can be naturally identified with $\m$, by
\[
T(x,y)=-x\cdot y,\qquad R(x,y)z=-[x,y,z],
\]
for any $x,y,z\in \m$ (see \cite[Theorem 10.3]{Nom54}).

Moreover, Nomizu showed too that the affine connections on manifolds
with parallel torsion and curvature are locally equivalent to
canonical connections  on reductive homogeneous spaces.
\textbf{}
Following the main ideas and results in \cite{Nom54}, Yamaguti introduced in \cite{Yam58}
what he called the \emph{general Lie triple systems}, later
renamed as \emph{Lie triple algebras} in \cite{Kik75}. We will
follow here the notation in \cite[Definition 5.1]{KinWei}, and will
call these systems \emph{Lie-Yamaguti algebras}:

\begin{definition}\label{df:LY}
A \emph{Lie-Yamaguti algebra} $(\m,x\cdot y,[x\,,y\,,z\,])$  ({\em
LY-algebra} for short) is a vector space $\m$ equipped with a
bilinear operation $\cdot : \m\times\m\rightarrow \m$ and a
trilinear operation $[\,,\,,\,]: \m\times\m\times\m\rightarrow \m$
such that, for all $x,y,z,u,v,w\in \m$:
\begin{enumerate}
\item[(LY1)] $x\cdot x=0$,
\item[(LY2)] $[x,x,y]=0$,
\item[(LY3)] $\sum_{(x,y,z)}\Bigl([x,y,z]+(x\cdot y)\cdot
z\Bigr)=0$,
\item[(LY4)] $\sum_{(x,y,z)}[x\cdot y,z,t]=0$,
\item[(LY5)] $[x, y,u\cdot v]=[x,y,u]\cdot v+u\cdot [x,y,v]$,
\item[(LY6)]
$[x,y,[u,v,w]]=[[x,y,u],v,w]+[u,[x,y,v],w]+[u,v,[x,y,w]]$.
\end{enumerate}
\end{definition}

\noindent Here $\sum_{(x,y,z)}$ means the cyclic sum on $x,y,z$.
\smallskip

The LY-algebras with $x\cdot y=0$ for any $x,y$ are exactly the Lie
triple systems, closely related with symmetric spaces, while the
LY-algebras with $[x,y,z]=0$ are the Lie algebras.  Less known
examples  can be found in \cite{BDE} where a detailed analysis on
the algebraic structure of LY-algebras arising from homogeneous
spaces which are quotients of the compact Lie group $G_2$ is given.
\smallskip

These nonassociative binary-ternary algebras have been treated by
several authors in connection with geometric problems on homogeneous
spaces \cite{Kik79,Kik81,Sag65,Sag68,SagWin}, but no much
information on their algebraic structure is available yet.
\smallskip

Following \cite{BEM09}, given a Lie-Yamaguti algebra $(\m,x\cdot y,[x,y,z])$ and any two
elements $x,y\in\m$, the linear map $D(x,y):\m\rightarrow\m$,
$z\mapsto D(x,y)(z)=[x,y,z]$ is, due to (LY5) and (LY6), a
derivation of both the binary and ternary products. These
derivations will be called \emph{inner derivations}. Moreover, let
$D(\m,\m)$ denote the linear span of the inner derivations. Then
$D(\m,\m)$ is closed under commutation thanks to (LY6). Consider the
vector space $\g(\m)=D(\m,\m)\oplus\m$, and endow it  with the
anticommutative multiplication given, for any $x,y,z,t\in \m$, by:
\begin{equation}\label{eq:gm}
\begin{split}
&[D(x,y),D(z,t)]= D([x,y,z],t)+D(z,[x,y,t]),\\
&[D(x,y),z]=D(x,y)(z)=[x,y,z],\\
&[z,t]=D(z,t)+z\cdot t.
\end{split}
\end{equation}
Note that the Lie algebra $D(\m,\m)$ becomes a subalgebra of
$\g(\m)$.

Then it is straightforward \cite{Yam58} to check that $\g(\m)$ is a
Lie algebra, called the \emph{standard enveloping Lie algebra} of
the Lie-Yamaguti algebra $\m$. The binary and ternary products in
$\m$ coincide with those given by \eqref{eq:binter}, where
$\h=D(\m,\m)$.

Given a Lie algebra $\g$ and a subalgebra $\h$, the pair $(\g,\h)$
will be said to be a \emph{reductive pair} (see \cite{Sag68}) if
there is a complementary subspace $\m$ of $\h$ with
$[\h,\m]\subseteq \m$. The decomposition $\g=\h\oplus\m$ will then
be called a \emph{reductive decomposition} of the Lie algebra $\g$
and  \emph{symmetric (reductive) decomposition} if the additional
condition $[\m,\m]\subseteq \h$ holds. In the latter case, we shall
refer to the pair $(\g,\h)$ as a {\em symmetric (reductive) pair}.
In particular, given a LY-algebra $(\m, x\cdot y,[x,y,z])$, the pair
$\bigl(\g(\m),D(\m,\m)\bigr)$ is a reductive pair and the pair is
symmetric in case $x\cdot y=0$.
\smallskip

As mentioned above, the Lie triple systems are precisely those
LY-algebras with trivial binary product. So they are related to symmetric decompositions and correspond to the
symmetric homogeneous spaces. Following \cite[\S 16]{Nom54}, a
symmetric homogeneous space $G/H$ is said to be \emph{irreducible}
if the action of $\ad\h$ on $\m$ is irreducible, where
$\g=\h\oplus\m$ is the canonical decomposition of the Lie algebra
$\g$ of $G$. This suggests the following definition (see \cite[Definition 1.2]{BEM09}):

\begin{definition}\label{df:irreducible}
A Lie-Yamaguti algebra $(\m, x\cdot y,[x,y,z])$ is said to be
\emph{irreducible} if $\m$ is an irreducible module for its Lie
algebra of inner derivations $D(\m,\m)$.
\end{definition}

The irreducible Lie-Yamaguti algebras constitute the algebraic
counterpart to the isotropy irreducible homogeneous spaces
considered in \cite{Wolf}. Concerning these irreducible LY-algebras
over algebraically closed fields of characteristic zero, it is not
difficult to prove (see \cite[Proposition 1.3, Theorem 2.1]{BEM09})
the following basic structure results:

\begin{theorem}\label{th:estructura}
Let $(\m, x\cdot y,[x,y,z])$ be an irreducible  LY-algebra. Then
$D(\m,\m)$ is a semisimple and maximal subalgebra of the standard
enveloping Lie algebra $\g(\m)$. Moreover, $\g(\m)$ is simple in
case $\m$ and $D(\m,\m)$ are not isomorphic as $D(\m,\m)$-modules. \qed
\end{theorem}

\begin{proposition}\label{pr:envueltasimple}
Let $\g=\h\oplus \m$ be a reductive decomposition of a simple Lie
algebra $\g$, with $\m\ne 0$. Then $\g$ and $\h$ are isomorphic,
respectively, to the standard enveloping Lie algebra and the inner
derivation algebra of the Lie-Yamaguti algebra $(\m,x\cdot
y,[x,y,z])$ given by \eqref{eq:binter}. Moreover, in case $\h$ is
semisimple and $\m$ is irreducible as a module for $\h$, either $\h$
and $\m$ are isomorphic as $\ad \h$-modules or $\m=\h^\perp$, the
orthogonal complement of $\h$ relative to the Killing form of $\g$.  \qed
\end{proposition}

From Theorem \ref{th:estructura}, in \cite[Section 2]{BEM09} it is
proved that the classification of irreducible LY-algebras splits
into three non overlapping types:
\begin{equation}\label{tipos}
\begin{array}{ll}

\textsc{Adjoint Type:}&\textrm{$\m$ is the adjoint module for $D(\m,\m)$,}\\

\textsc{Non-Simple Type:} &\textrm{$D(\m,\m)$ is not simple,}\\

\textsc{Generic Type:}&\textrm{Both $\g(\m)$ and $D(\m,\m)$ are
simple.}
\end{array}
\end{equation}

The LY-algebras of Adjoint Type are just the simple Lie algebras
(see \cite[Theorem 2.4]{BEM09}) and those of Non-Simple Type can be
described through reductive decompositions modeled by a Generalized
Tits Construction from \cite{BenZel} using quaternions, octonions
and simple Jordan algebras as basic ingredients (see \cite[Theorems
4.1, 4.4]{BEM09}).

In the Generic Type, $\m$ and $D(\m,\m)$ are not isomorphic as
$\ad_{\g(\m)} D(\m,\m)$-modules, so following Proposition
\ref{pr:envueltasimple}, the classification of the irreducible
LY-algebras of this type is equivalent to the determination of the
reductive decompositions $\g=\h\oplus \m$ satisfying the following
conditions:
\begin{equation}\label{condicionescasonosimpleI}
\begin{array}{rl}
&  \mathrm{(a)}\quad
\textrm{$\g$ is a simple Lie algebra,}\\
&\textrm{(b)}\quad
\textrm{$\h$ is a simple subalgebra of $\g$,} \\
& \mathrm{(c)}\quad \textrm{$\m$ is an irreducible $\ad\h$-module
(in particular $\m \neq 0$).}
\end{array}
\end{equation}
Note that the previous conditions imply
\begin{equation}\label{condicionescasonosimpleII}
\begin{array}{rl}
&  \mathrm{(d)}\quad
\textrm{$\m$ is $\h^\perp$ (orthogonal with respect to the Killing form of $\g$),}\\
&\textrm{(e)}\quad
\textrm{$\h$ is a maximal subalgebra of $\g$,} \\
& \mathrm{(f)}\quad \textrm{$\m$ is not the adjoint module.}
\end{array}
\end{equation}

The purpose in this paper is the classification of the LY-algebras
of Generic Type while, at the same time, their close connections to
some well-known nonassociative algebraic systems will be
highlighted. It will be shown that most part of irreducible
LY-algebras of this type appear inside simple Lie algebras as
orthogonal complements of subalgebras of derivations of Lie and
Jordan algebras and triple systems, Freudenthal and orthogonal
triple systems or Jordan and anti-Jordan pairs.

The paper is structured as follows. Section 2 is devoted to
determine the irreducible LY-algebras inside reductive
decompositions of simple special linear Lie algebras (classical
Cartan type $A_n$). The classification of these LY-algebras flows
parallel to the classification of the simple Jordan linear pairs and
the so called anti-Jordan pairs. Following a similar pattern,
Sections 3 and 4 provide LY-algebras appearing inside orthogonal and
symplectic simple Lie algebras (Cartan types $B_n, D_n$ and $C_n$)
through the classification of simple Lie triple systems, orthogonal
and symplectic triple systems. Irreducible LY-algebras inside
exceptional Lie algebras of types $G_2, F_4, E_6, E_7$ and $E_8$ are
the goal of Section 5. In this case, the classification can be
transferred from the complex field. Section 6 is an appendix section
where definitions and classifications of the different pairs and
triple systems related to irreducible LY-algebras are included. The
paper ends with an epilogue section that summarizes the
classification results obtained in this paper and in the previous
one \cite{BEM09}.

Throughout this paper, \emph{all the algebraic systems considered
will be assumed to be finite dimensional over an algebraically
closed ground field $k$ of characteristic zero}. The symbol $\oplus$
denotes the direct sum of subspaces and $\otimes$ the tensor product
of $k$-subspaces unless otherwise stated. Basic notation and
terminology on representation theory of Lie algebras follows
\cite{Hum72}.

\section{Special linear case}
As mentioned in the Introduction, the irreducible LY-algebras of
Generic Type appear as orthogonal complements of maximal simple
subalgebras of simple Lie algebras. In this section we classify
these  systems in case their standard enveloping algebras are
(simple) special linear Lie algebras  $\spl(V)$ (or $\spl_{\dimens
V}(k)$ if a basis for $V$ is fixed).

Any reductive decomposition $\spl(V)=\h\oplus \m$ satisfying (a),
(b) and (c) in (\ref{condicionescasonosimpleI}) presents two
elementary restrictions:
\begin{itemize}
\item[$\bullet$] $\dimens\, V \geq 3$ (the smallest simple Lie algebra
is the three dimensional algebra $\spl_2(k)$).
\item[$\bullet$] $V$ is irreducible as a module for $\h$. Otherwise
we would have $V=V_1\oplus V_2$ and $\h\subseteq \{x\in\spl(V) :
f(V_i)\subseteq V_i, \tr f\mid_{V_i}=0\}$, but this subspace is
properly contained in the subalgebra $\{f\in\spl(V) : f(V_i)\subseteq V_i\}$.
This is not possible by the maximality of the subalgebra $\h$, following
condition (e) from (\ref{condicionescasonosimpleII}).
\end{itemize}

The previous restrictions allow us to introduce trilinear products involving $V$ and its dual $\h$-module $V^*$, in such a way that the pair $(V,V^*)$ is endowed with a (linear) Jordan or anti-Jordan pair structure (see \cite{Loos} and \cite{FauFe80} for definitions or Subsection \ref{pares} of this paper) and $\h$ can be viewed as the derived subalgebra of the inner derivation algebra of the induced pair.

Since $\spl(V)$ is embedded in $\gl(V)$, and this is a module for its subalgebra $\h$, we can consider the standard isomorphism of $\h$-modules
\begin{equation}\label{isogl(v)}
V\otimes V^*\cong \gl(V)=\h\oplus \m \oplus kI_V
\end{equation}
given by $x\otimes \varphi \mapsto \varphi(-)x$. ($I_V$ denotes the identity map of the vector space $V$.)

Now let $d_{x, \varphi}$ be the projection of $\varphi(-)x$ onto
$\h$ and, for a fixed $\xi\in k$, let us define the $\h$-invariant
triple products

\begin{equation}\label{triple+-sl}
\begin{array}{rcl}
V\otimes V^*\otimes V & \to  &V\\
x\otimes\varphi\otimes y & \mapsto &  \{x\varphi y\}_{\xi}:=d_{x,\varphi}(y)-\xi \varphi(x)y
\end{array}
\end{equation}
\begin{equation}\label{triple--sl}
\begin{array}{rcl}
V^*\otimes V\otimes V^* & \to  &V^*\\
\varphi\otimes x\otimes\psi & \mapsto & \{\varphi x \psi\}_{\xi} :=  
\psi\circ d_{x,\varphi}-\xi \varphi(x)\psi\\
\end{array}
\end{equation}
\medskip

Products (\ref{triple+-sl}) and (\ref{triple--sl}) are related, for all $\varphi,\psi \in V^*$,  $x,y \in V$, by
\begin{equation}\label{triple-producto-relacion}
\{\varphi x \psi \}_\xi(y)=\psi \circ \{x \varphi y\}_\xi
\end{equation}
and the subalgebra $\h$ can be described as
\begin{equation}\label{triple-h}
\h=\eespan\langle d_{x, \varphi}: x \in V, \varphi \in V^* \rangle
\end{equation}
Then, we have the following result:

\begin{lemma}\label{casoparJordan}
For a given reductive decomposition  $\spl(V)=\h\oplus \m$ which satisfies {\em (a), (b), (c)} in {\em (\ref{condicionescasonosimpleI})}, consider the vector spaces $U^+=V$ and $U^-=V^*$, there exists a nonzero scalar $\xi\in k$ such that the pair $\mathcal{U}=(U^+,U^-)$ is either a simple Jordan pair under the triple products $\{x_{\sigma}y_{-\sigma}z_{\sigma}\}_{\xi}$ defined in  {\em (\ref{triple+-sl})} and {\em (\ref{triple--sl})} for {\em $\sigma=\pm$}, or a contragredient simple anti-Jordan pair with $\langle x_{\sigma}y_{-\sigma}z_{\sigma}\rangle:= \sigma \{x_{\sigma}y_{-\sigma}z_{\sigma}\}_{\xi}$ as triple products. Moreover $\h$ is the linear subalgebra
\begin{equation}\label{las-h-sl}
\h=\eespan\langle \{x\varphi \,\cdot  \}_{\xi}-\xi\varphi(x)I_V: x\in V, \varphi \in V^*\rangle
\end{equation}
which, up to isomorphism, turns out to be the derived subalgebra of the inner derivation Lie algebra of the corresponding pair and $\m$ is $\h^\perp$, the orthogonal complement of $\h$ with respect to the Killing form of $\spl(V)$.
\end{lemma}
\begin{proof}

First we shall check that, for an arbitrary $\xi$, the $\xi$-products in (\ref{triple+-sl}) and (\ref{triple--sl}) satisfy the identity
\begin{equation}\label{Jp2-sl}
\begin{array}{l}
\{x_{\sigma}y_{-\sigma}\{u_{\sigma}v_{-\sigma}w_{\sigma}\}_{\xi}\}_{\xi}
=\{\{x_{\sigma}y_{-\sigma}u_{\sigma}\}_{\xi}v_{-\sigma}w_{\sigma}\}_{\xi}\\
 \qquad\qquad\qquad-\{u_{\sigma}\{y_{-\sigma}x_{\sigma}v_{-\sigma}\}_{\xi}w_{\sigma}\}_{\xi}
  +\{u_{\sigma}v_{-\sigma}\{x_{\sigma}y_{-\sigma}w_{\sigma}\}_{\xi}\}_{\xi}
\end{array}
\end{equation}
For $x=x_+$, $u=u_+\in U^+=V$ and $\varphi=y_-$, $\psi=v_-\in U^-=V^*$, the map
\begin{equation}\label{Lx-fi}
L:V\otimes V^*\to \End(V)
\end{equation}
defined by $x \otimes \varphi \mapsto L_{x,\varphi}=\{x\varphi-\}_\xi$ is an ($\h\oplus k1$)-module homomorphism. From (\ref{triple-producto-relacion}) we have $\{\varphi x \psi \}_\xi=\psi \circ L_{x,\varphi}$ that easily yields to
\begin{equation}\label{equiv}
[L_{x,\varphi},L_{u,\psi}]=L_{L_{x,\varphi}(u),\psi}-L_{u, \psi \circ L_{x,\varphi}}
\end{equation}
But (\ref{equiv}) is equivalent to (\ref{Jp2-sl}) for $\sigma=+$.

In case $\varphi=x_-$, $\psi=u_-\in U^-=V^*$ and $x=y_+$, $y=v_+\in U^*=V$, identity (\ref{Jp2-sl}) follows from the ($\h\oplus k1$)-module homomorphism
\begin{equation}\label{Lx-fi2}
\hat{L}:V^*\otimes V\to \End(V^*)
\end{equation}
given by $\varphi \otimes x \mapsto \hat{L}_{\varphi,x}=\{\varphi x-\}_\xi$.

On the other hand, any product defined as in (\ref{triple+-sl}) is
in $\Hom_{\h}(V \otimes V^*\otimes V,V)$, so we must look at the
previous subspace in order to get our result. Since $\h$, $\m$ and
$kI_V$ are irreducible and non-isomorphic $\h$-modules,
$$
\begin{array}{rcl}
\Hom_{\h}(V \otimes V^*\otimes V,V) & \cong & \Hom_{\h}(V \otimes V^*,V \otimes V^*)\\
& \cong & \Hom_{\h}(\h,\h) \oplus \Hom_{\h}(\m,\m) \oplus \Hom_{\h}(kI_V,kI_V)
\end{array}
$$
Then $\Hom_{\h}(V \otimes V^*\otimes V,V)$ is a $3$-dimensional vector space from Schur's Lemma. Now, using the alternative decomposition
$$
\begin{array}{rcl}
\Hom_{\h}(V \otimes V^*\otimes V,V) &\cong & \Hom_{\h}(V \otimes V,V\otimes V)\\
& \cong & \Hom_{\h}(S^2 V,V\otimes V)\oplus \Hom_{\h}(\wedge^2 V,V\otimes V)
\end{array}
$$
we get that either $\Hom_{\h}(\wedge^2 V,V\otimes V)$ or  $\Hom_{\h}(S^2 V,V\otimes V)$ is a one-dimensional subspace ($S^2V$ and $\wedge^2V$ stand for the second symmetric and alternating power of $V$ respectively). Hence, two different situations appear:
\medskip

a) $\Hom_{\h}(\wedge^2 V,V\otimes V)$ is one-dimensional

\noindent In this case, the set $\Hom_{\h}(\wedge^2 V\otimes V^*,V)$ is also one-dimensional and, since $\dimens V \geq 3$, we can take the nonzero map
$$
x \otimes y \otimes \varphi  - y \otimes x \otimes \varphi \mapsto \varphi(x)y- \varphi(y)x
$$
as generator. So
$$
d_{x,\varphi}(y)- d_{y,\varphi}(x)=\xi (\varphi(x)y- \varphi(y)x)
$$
for some $\xi \in F$ and therefore, for $x,y\in V$, $\varphi\in V^*$ we have the identity
\begin{equation}\label{eq3.6}
\{x\varphi y\}_{\xi}= \{y\varphi x\}_{\xi}
\end{equation}
On the other hand, since $\Hom_{\h}(\wedge^2 V,V\otimes V)$ is assumed to be one-dimensional and $V\otimes V=S^2 V\oplus \wedge^2 V$, we have that $\Hom_{\h}(\wedge^2 V,S^2 V)=\Hom_{\h}(S^2 V,\wedge^2 V)=\Hom_{\h}(S^2 V\otimes \wedge^2 V^*,k)=0$. This latter condition implies that the restriction of the map $x\otimes y \otimes \varphi \otimes \psi \mapsto \psi(\{x \varphi y\}_{\xi})$ on $S^2 V\otimes \wedge^2 V^*$ is zero. So, using (\ref{triple-producto-relacion}) and taking into account that the base field is of characteristic zero, we get
$$
0=\psi(\{x\varphi y\}_{\xi})-\varphi(\{x\psi y\}_{\xi})=(\{\varphi x \psi\}_{\xi}-\{\psi x \varphi\}_{\xi})(y)
$$
which implies
\begin{equation}\label{eq3.6'}
\{\varphi x \psi\}_{\xi}=  \{\psi x \varphi\}_{\xi}
\end{equation}
Then, from (\ref{Jp2-sl}), (\ref{eq3.6}) and (\ref{eq3.6'}), we obtain that $(\mathcal{U},\{x_{\sigma}y_{-\sigma}z_{\sigma}\}_{\xi})$ is a Jordan pair.

b) $\Hom_{\h}(S^2 V,V\otimes V)$ is one-dimensional

\noindent In this case, analogous arguments but for symmetric powers, give us
$$
d_{x,\varphi}(y)+ d_{y,\varphi}(x)=\xi (\varphi(x)y+ \varphi(y)x)
$$
for some $\xi\in k$, which yields to

\begin{equation}\label{eq3.6b}
\{x\varphi y\}_{\xi}+ \{y\varphi x\}_{\xi}=0
\end{equation}
and
\begin{equation}\label{eq3.6'b}
\{\varphi x \psi\}_{\xi}+  \{\psi x \varphi\}_{\xi}=0
\end{equation}
for $x,y\in V$, $\varphi, \psi \in V^*$. Now for $\sigma=\pm$, the products $\langle x_{\sigma}y_{-\sigma}z_{\sigma}\rangle=\sigma \{x_{\sigma}y_{-\sigma}z_{\sigma}\}_{\xi}$ satisfy
\begin{equation}\label{aJp1-sl}
\langle x_{\sigma}y_{-\sigma}z_{\sigma}\rangle=-\langle z_{\sigma}y_{-\sigma}x_{\sigma}\rangle
\end{equation}
and using (\ref{Jp2-sl}):

\begin{equation}\label{aJp2-sl}
\begin{array}{l}
\langle x_{\sigma}y_{-\sigma} \langle u_{\sigma}v_{-\sigma}z_{\sigma} \rangle \rangle
 = \sigma^2\{x_{\sigma}y_{-\sigma}\{ u_{\sigma}v_{-\sigma}z_{\sigma}\}_{\xi}\}_{\xi}\\
  \qquad  \qquad = \sigma^2\{\{x_{\sigma}y_{-\sigma}u_{\sigma}\}_{\xi}v_{-\sigma}z_{\sigma}\}_{\xi}-
\sigma^2\{u_{\sigma}\{y_{-\sigma}x_{\sigma}v_{-\sigma}\}_{\xi}z_{\sigma}\}_{\xi}\\
\qquad  \qquad \quad
+\sigma^2\{u_{\sigma}v_{-\sigma}\{x_{\sigma}y_{-\sigma}z_{\sigma}\}_{\xi}\}_{\xi}\\
\qquad \qquad = \langle \langle x_{\sigma}y_{-\sigma}u_{\sigma} \rangle v_{-\sigma}z_{\sigma} \rangle\\
\qquad \qquad\quad +\langle u_{\sigma} \langle
y_{-\sigma}x_{\sigma}v_{-\sigma} \rangle z_{\sigma} \rangle +
\langle u_{\sigma}v_{-\sigma} \langle
x_{\sigma}y_{-\sigma}z_{\sigma} \rangle \rangle
\end{array}
\end{equation}

Now identities (\ref{aJp1-sl}) and (\ref{aJp2-sl}) prove that
$(\mathcal{U}, \langle x_{\sigma}y_{-\sigma}z_{\sigma}\rangle)$ is an anti-Jordan pair.

Following \cite[Section 1]{FauFe80}, for a given Jordan or anti-Jordan pair with triple products $a_\sigma b_{-\sigma} c_\sigma=D_\sigma(a_\sigma,b_{-\sigma})(c_\sigma)$, the so called inner derivation algebra is the Lie algebra spanned by the (inner) derivations $D(a_+,b_{-})=(D_+(a_+,b_-),-D_-(b_-,a_+))$ in the Jordan pair case or $D(a_+,b_-)=(D_+(a_+,b_-)$, $D_-(b_-,a_+))$ in the anti-Jordan pair case. In this way, comparing inner derivation maps for the pair and anti-pair obtained from a) or b), we arrive at the relationship
$$
(\langle x \varphi -\rangle,\langle \varphi  x -\rangle)=(\{x \varphi -\}_\xi,-\{\varphi x -\}_\xi)
$$
Thus in both cases,
$$
\Inder\, \mathcal{U}= \eespan\langle (\{x\varphi -\}_{\xi},-\{\varphi x -\}_{\xi}): x \in V , \varphi \in V^*\rangle
$$
Now, the Lie algebra $\Inder\, \mathcal{U} $ is isomorphic to $\h$ in case $\xi =0$ or $\h\oplus kI_V$ otherwise. The previous assertion follows from the map $\h\oplus kI_V\to \Inder\  \mathcal{U}$, given by $d \mapsto (d,-\tilde{d})$, where $\tilde{d}(\varphi)=\varphi\circ d$. Moreover, since $V$ and $V^*$ are $\h$-irreducible, according to \cite[Proposition 1.2]{FauFe80}, $\mathcal{U}$ is a simple Jordan pair or anti-Jordan pair. Attached to the anti-Jordan pair structure, we have the nondegenerate bilinear map $V\otimes V^* \to k$ defined by $x\otimes \varphi \mapsto \varphi(x)$ which, because of  (\ref{triple-producto-relacion}) satisfies
$$
\psi(\langle x \varphi y\rangle)+\langle\varphi x\psi\rangle(y)=
\psi(\{x \varphi y\}_{\xi})-\{\varphi x\psi\}_{\xi}(y)=0
$$
$x,y\in V$, $\varphi,\psi\in V^*$. Then the anti-Jordan pair is contragredient following \cite[Section 2]{FauFe80}.  Now, from \cite[Theorem 2]{Mey84} and \cite[Sections 2 and 3]{FauFe80}, the inner derivation Lie algebra of either a simple Jordan pair or a simple contragredient anti-Jordan pair is never simple (see Tables \ref{simple-jordan} and \ref{simple-antijordan} in the Appendix section for a complete description of these algebras). So, up to isomorphisms, the Lie algebra $\Inder\, \mathcal{U}$ is $\h \oplus kI_V$ which proves that $\xi\neq 0$. So the derived subalgebra $\Inder_{0}\, \mathcal{U}=[\Inder\, \mathcal{U},\Inder\, \mathcal{U}]$ is isomorphic to $\h$. Moreover, using (\ref{triple-h}), the algebra $\h$ is spanned by the zero trace maps $d_{x,\varphi}=\{x\varphi \, \cdot\}_{\xi}+\xi\varphi(x)I_V$, therefore  $\xi\varphi(x)=-\frac{\tr(\{x\varphi \, \cdot\}_{\xi})}{\dim V}$. The final assertion on $\m$ follows from condition (d) in (\ref{condicionescasonosimpleII}).
\end{proof}

Now, we can establish the main result for the generic $\spl$-case:

\begin{theorem}\label{descripcion-sl}
Let $(\m, a\cdot b, [a,b,c])$ be an irreducible LY-algebra of generic type and standard enveloping Lie algebra of type $\spl$. Then either:
\begin{itemize}
\item[(i)] There is a vector space $V$ and an involution on the associative algebra $\End(V)$ such that $\m$ is, up to isomorphism, the simple Lie triple system consisting of the zero trace symmetric elements in $\End(V)$, with the natural triple product $[a,b,c]=[[a,b],c]$ (inside $\spl(V)$). Moreover, $\dim V\geq 5$ if the involution is orthogonal, and $\dim V\geq 4$ if it is symplectic. In particular, the binary product $a\cdot b$ is trivial.

\item[(ii)] There is a simple Jordan triple system $J$ of one of the following types:
    \begin{enumerate}
    \item the subspace of $n\times n$ symmetric matrices for $n\geq 2$ with the triple product $\{xyz\}=xy^tz+zy^tx$,
    \item the subspace of $n\times n$ skew-symmetric matrices for $n\geq 5$ again with the triple product $\{xyz\}=xy^tz+zy^tx$,
    \item the subspace of $1\times 2$-matrices over the algebra of octonions $\cO$ with the triple product  $\{xyz\}=x({\bar y}^tz)+z({\bar y}^tx)$,
    \item the exceptional Jordan algebra $\cH_3(\cO)$ (multiplication denoted by juxtaposition) with its triple product $\{xyz\}=x(zy)+z(xy)-(zx)y$.
    \end{enumerate}
    such that, up to isomorphism, $\g(\m)=\spl(J)$, $\h=D(\m,\m)=\mathcal{L}_{_0}(J)=$ $[\cL(J),\cL(J)]$, where $\cL(J)=\eespan\langle \{xy.\}; x,y\in J\rangle$. Here the LY-algebra $\m$ appears as the orthogonal complement to $\h$ in $\g(\m)$ relative to the Killing form, with the binary and ternary products in \eqref{eq:binter}.
\end{itemize}

There are no isomorphisms among the LY-algebras in the different items above.
\end{theorem}

\begin{proof}
According to Lemma \ref{casoparJordan},  $\m=\h^\perp$ where $\h$ is
as described in (\ref{las-h-sl}) and $V=\mathcal{U}^+$, the
$(+)$-component of either a suitable  simple Jordan pair or  a
simple contragredient anti-Jordan pair
$\mathcal{U}=(\mathcal{U}^+,\mathcal{U}^-)$ with triple products
$\{x_\sigma y_\sigma z_\sigma\}$ for $\sigma=\pm$. We also note that
the subalgebra $\h$ is described by means of the $(+)$-product
operators $D_+(x_+,y_{-})=\{x_{+} y_{-}\,  \cdot\, \}$, $x_+ \in
\mathcal{U}^+$ and $y_-\in \mathcal{U}^-$. Moreover, for
contragredient anti-Jordan pairs we must have the isomorphism
$(\mathcal{U}^+)^*\cong \mathcal{U}^-$, as modules for $\Inder\,
\mathcal{U}$. (Any simple Jordan pair is contragredient, but the
assertion is not true for anti-Jordan pairs as shown in
\cite[Example 2.7]{FauFe80}.)

Conditions (b) and (c) in (\ref{condicionescasonosimpleI}) and the
previous  initial restrictions on irreducibility and dimension of
$V$ imply that we must look only for simple Jordan or simple
contragredient anti-Jordan pairs
$\mathcal{U}=(\mathcal{U}^+,\mathcal{U}^-)$ such that:
\begin{itemize}
\item[(1)] the derived subalgebra
$\Inder_0\, \mathcal{U}=[\Inder\, \mathcal{U},\Inder\, \mathcal{U}]$
of the inner derivation algebra is simple,
\item[(2)] $\mathcal{U}^+=V(\lambda)$ is an irreducible module for the $(+)$-component
of $\Inder_0\, \mathcal{U}$ and $\dimens \mathcal{U}^+ \geq 3$,
\item[(3)] $\mathcal{U}^+\otimes \mathcal{U}^-=V(\lambda)\otimes V(\lambda^*)$
has just two non-trivial irreducible components which correspond to
the module decomposition for $\spl(\mathcal{U}^+)$.
\end{itemize}

Isomorphisms between either Jordan pairs or anti-Jordan pairs provide isomorphic LY-algebras. A look at the classification of simple Jordan pairs in \cite[Theorem 17.12]{Loos} and of simple anti-Jordan pairs in \cite[Section 3 and 4]{FauFe80} (both classifications are outlined in the Appendix --subsection \ref{pares} and Table \ref{simple-jordan}--) show that the possibilities for Jordan pairs $\mathcal{U}=(\mathcal{U}^+,\mathcal{U}^-)$ satisfying conditions (1) and (2) above are the following:
\begin{description}
\settowidth{\labelwidth}{XX}%
\setlength{\leftmargin}{30pt}
\item[$(\cM_{p,1}(k), \cM_{p,1}(k))_{p \geq 3}$]\null\quad\newline
    In this case, the pair $(\Inder_0\, \mathcal{U},\mathcal{U}^+)$ is, up to isomorphism, the pair $\bigl(\spl_p(k),V(\lambda_1)\bigr)$ (recall that we follow the notations of \cite{Hum72}). Then, as modules for $\spl_p(k)$, $\cU^+\otimes\cU^- = V(\lambda_1)\otimes V(\lambda_{p-1})=V(\lambda_1+\lambda_{p-1})\oplus V(0)$.\null\quad\newline
    Therefore, this case must be discarded (see \eqref{isogl(v)}).

\item[$(\cA_n(k),\cA_n(k))_{n \geq 5}$]\null\quad\newline
    Here the pair $(\Inder_0\, \mathcal{U},\mathcal{U}^+)$ is, up to isomorphism, $(\spl_n(k), V(\lambda_2))$ and, as modules for $\spl_n(k)$, $\cU^+\otimes\cU^- =V(\lambda_2)\otimes V(\lambda_{n-2})=V(\lambda_1+\lambda_{n-1})\oplus V(\lambda_2+\lambda_{n-2})\oplus V(0)$.

\item[$(\cH_n(k),\cH_n(k))_{n \geq 2}$]\null\quad\newline
    Here the pair $(\Inder_0\, \mathcal{U},\mathcal{U}^+)$ is, up to isomorphism, $(\spl_n(k),V(2\lambda_1))$   and, as modules for $\spl_n(k)$, $\cU^+\otimes\cU^- =   V(2\lambda_1)\otimes V(2\lambda_{n-1})=V(2\lambda_1+2\lambda_{n-1})\oplus V(\lambda_1+\lambda_{n-1})\oplus V(0)$.

\item[$(k^n,k^n)_{n\geq 5}$]\null\quad\newline
    Here the pair $(\Inder_0\, \mathcal{U},\mathcal{U}^+)$ is, up to isomorphism, $(\so_n(k),V(\lambda_1))$   and, as modules for $\so_n(k)$, $\cU^+\otimes\cU^- = V(\lambda_1)\otimes V(\lambda_1)$ is equal to $V(2\lambda_1)\oplus V(\lambda_2)\oplus V(0)$ for $n \geq 7$ and for $n=5, 6$ it decomposes as $V(2\lambda_1)\oplus V(2\lambda_2)\oplus V(0)$ and  $V(2\lambda_1)\oplus V(\lambda_2+\lambda_3)\oplus V(0)$ respectively.

\item[$(\cM_{1,2}(\cO),\cM_{1,2}(\cO))$]\null\quad\newline
    Here the pair $(\Inder_0\, \mathcal{U},\mathcal{U}^+)$ is, up to isomorphism, $(\so_{10}(k),V(\lambda_4))$   and, as modules for $\so_{10}(k)$, $\cU^+\otimes\cU^- =   V(\lambda_4)\otimes V(\lambda_5)=V(\lambda_4+\lambda_5)\oplus V(\lambda_2)\oplus V(0)$.

\item[$(\cH_3(\cO),\cH_3(\cO))$]\null\quad\newline
    Here the pair $(\Inder_0\, \mathcal{U},\mathcal{U}^+)$ is, up to isomorphism, $(E_6,V(\lambda_1))$   and, as modules for $E_6$, $\cU^+\otimes\cU^- =   V(\lambda_1)\otimes V(\lambda_6)=V(\lambda_1+\lambda_6)\oplus V(\lambda_2)\oplus V(0)$.
\end{description}

For anti-Jordan pairs, the results in Subsection \ref{pares} and Table \ref{simple-antijordan} show that the series $\mathcal{U}=(\cM_{p,1}(k), \cM_{p,1}(k))$ for ${p\ge 3}$ with $(\spl_p(k), V(\lambda_1))$ and $\mathcal{U}=(k^{2n}, k^{2n})$ with $(\frsp_{2n}(k), V(\lambda_1))$ for ${n\ge 2}$ are the unique possibilities. The decomposition $\mathcal{U}^+\otimes \mathcal{U}^-$ as $\spl_p(k)$-module in the first case is analogous to the corresponding series of Jordan pairs, and hence this case must be discarded. For the anti-Jordan pairs $\mathcal{U}=(k^{2n}, k^{2n})$, the decomposition as modules for $\frsp_{2n}(k)$, is given by
\begin{description}
\settowidth{\labelwidth}{XX}%
\setlength{\leftmargin}{30pt}
\item[$(k^{2n}, k^{2n})_{n \geq 2}$]\null\quad\newline
    $V(\lambda_1)\otimes V(\lambda_1)=V(2\lambda_1)\oplus V(\lambda_2)\oplus V(0)$
\end{description}

\smallskip

The Jordan pairs and anti-Jordan pairs of type $\mathcal{U}=(k^n,k^n)$ present a special common feature. Both structures can be described as a pair $(V,V)$  where $V$ is a vector space endowed with a nondegenerate $\epsilon$-symmetric form $b$ and with triple products $\{xyz\}=b(x,y)z+b(y,z)x-\epsilon b(x,z)y$, where $\epsilon =1$ for Jordan pairs and $\epsilon=-1$ for anti-Jordan pairs (so  $\dimens V$ is even in the latter case). Moreover, the operators appearing in (\ref{las-h-sl}) are of the form
$$
d_{x,y}=\{x y \, \cdot \}-\frac{\tr(\{x y \, \cdot\})} {\dimens V}= b(y,z)x-\epsilon b(x,z)y
$$
and hence the subalgebra $\h=\eespan \langle d_{x,y}= b(y,\cdot)x-\epsilon b(x,\cdot)y: x,y \in V \rangle$ is the Lie algebra $\frso(V)$ in case $\epsilon=1$ and $\frsp(V)$ for $\epsilon=-1$. On the other hand, the map $f\mapsto f^*$, where $f^*$ is the adjoint map relative to the form $b$, induces an involution on the associative algebra $\End(V)$ for which $\h=\{f\in\End(V): f^*=-f\}=\cS(V,*)$ is just the Lie algebra $\frso(V)$ or $\frsp(V)$ and the set $\cJ=\cH(V,*)=\{f\in\End(V):f^*=f\}$, under the symmetrized product $f\cdot g=fg+gf$, is a central simple Jordan algebra. In this case, the decomposition $\gl(V)=\cS(V,*)\oplus\cH(V,*)$ is symmetric and its restriction to $\spl(V)$ provides the symmetric decomposition $\spl(V)=\h\oplus \cH(V,*)_0$, where $\cH(V,*)_0$ consist of the zero trace elements in $\cJ$. This symmetric decomposition satisfies (\ref{condicionescasonosimpleI}), so $\m=\h^\perp=\cH(V,*)_0$ and therefore the LY-algebra $\m$ has trivial binary product and the ternary  one is given by $[f,g,h]=[[f,g],h]=-(f,h,g)=-(f\cdot h)\cdot g+f\cdot(h\cdot g)=(g,f,h)$. This provides item (i) in the Theorem.


Finally, the remaining admissible Jordan pairs above are all of the form $(J,J)$ for a Jordan triple system $J$ and the Theorem follows.
\end{proof}


\section{Orthogonal case}

In this section we classify LY-algebras of Generic Type  in case
their standard enveloping is a (simple) orthogonal Lie algebra $\so(V,b)$, so $V$ is a vector space of dimension $\geq 5$ (note that $\so_4(k)$ is not simple and that $\so_3(k)$ is isomorphic to $\spl_2(k)$), endowed with a nondegenerate symmetric form $b$.

As in the previous section, we are looking for decompositions $\so(V,b)=\h\oplus\m$ in which conditions (a), (b) and (c) in
(\ref{condicionescasonosimpleI}) hold. Our discussion in the
$\so$-case will be based on the following elementary facts:
\begin{itemize}
\item[$\bullet$] Considering both $\so(V)$ and $V$ as modules for $\h$, the linear map $x\wedge y\mapsto \sigma_{x,y}=b(x,.)y-b(y,.)x$ defines an isomorphism of $\h$-modules:
\begin{equation}\label{sk-so}
\wedge^2V\cong\so(V,b)
\end{equation}
    from the second alternating power of $V$ onto the Lie algebra $\so(V)$.
\item[$\bullet$] Any tensor product of irreducible modules $V(\lambda)\otimes V(\mu)$ contains a (unique) copy of the irreducible module $V(\lambda +\mu)$. This copy is generated by $v=v_{\lambda}\otimes v_{\mu}$, the only vector (up to scalars) of (highest) weight $\lambda +  \mu$. (Here $v_\lambda$ denotes a nonzero vector of weight $\lambda$.) Moreover, in case $\lambda =\mu$ this copy is located inside the second symmetric power of $V(\lambda)$, that is:
\begin{equation}\label{copsim}
V(2\lambda)\subseteq S^2(V(\lambda)).
\end{equation}

\item[$\bullet$] For a given dominant weight $\lambda$ and any simple root $\alpha$ not orthogonal to $\lambda$ ($\langle \lambda,\alpha\rangle \neq 0$), the second alternating power $\wedge^2 V(\lambda)$ contains a (unique) copy of the irreducible module $V(2\lambda-\alpha)$. This copy is generated by $v=v_{\lambda}\otimes v_{\lambda-\alpha}-v_{\lambda -\alpha}\otimes v_{\lambda}$, the only vector (up to scalars) of (highest) weight $2\lambda -\alpha$. Hence,
\begin{equation}\label{copsk}
V(2\lambda-\alpha)\subseteq \wedge^2(V(\lambda)),\
\mathrm{in\  case}\  \langle \lambda, \alpha \rangle \neq 0
\end{equation}
\end{itemize}

\begin{lemma}\label{estructura-so}
For a given reductive decomposition $\so(V,b)=\h\oplus\m$ satisfying {\em
(a)}, {\em (b)} and {\em (c)} in {\em
(\ref{condicionescasonosimpleI})} with  $\dimens V \ge 5$, one has that, as a module for $\h$, either:
\begin{itemize}
\item[(i)] $V$ decomposes as $V=k v\oplus W$, an orthogonal sum of a trivial module $k v$ and an irreducible module $W$ with $\dimens W \ge 5$. In this case, the subalgebra $\h$ is $\h=\sigma_{W,W}=\eespan \langle \sigma_{x,y}:x,y\in W\rangle$, ($\sigma_{x,y}$ as in {\em (\ref{sk-so})}), so it is isomorphic to $\so(W,b)$ and for the subspace $\m$ we have that $\m=\sigma_{v,W}$. Moreover, the reductive decomposition
\begin{equation}\label{so-triples}
\so(kv\oplus W,b)=\sigma_{W,W}\oplus \sigma_{v,W}
\end{equation}
is symmetric; or

\item[(ii)] $V=V(m\lambda_i)$ is an irreducible module for $\h$ whose dominant weight is a multiple of the fundamental weight $\lambda_i$ relative to some system of simple roots $\Delta=\{\alpha_1, \dots, \alpha_n\}$, and one of the following holds:
\begin{itemize}
\item[(ii-a)] $\m=V(2m\lambda_i-\alpha_i)$ for some $i$, $1\leq i\leq n$.
\item[(ii-b)] $\h$ is a simple Lie algebra of type $B_3$, $V=V(\lambda_3)$ and $\m=V(\lambda_1)$
\item[(ii-c)] $\h$ is a simple Lie algebra of type $G_2$ and $V=\m=V(\lambda_1)$.
\end{itemize}

\end{itemize}

\end{lemma}

\begin{proof}
If $V$ is not irreducible as a module for $\h$, let $W$ be
a proper and irreducible $\h$-submodule. Assume first $b(W,W)\neq 0$, thus the restriction of $b$ to $W$ is nondegenerate by irreducibility of $W$, so $V$ decomposes as the orthogonal sum $V=W\oplus W^\perp$. Since $\h \subset \so(W)\oplus\so(W^\perp)\subset \so(V,b)$, the maximality of $\h$ (condition (e) in (\ref{condicionescasonosimpleII})) forces
$$
\h = \so(W,b)= \so(W)\oplus \so(W^\perp)
$$

But $\h$ is a simple Lie algebra, so $W^\perp $ must be one-dimensional. So we have the orthogonal decomposition
$V=kv\oplus W$ and we get the natural $\mathbb{Z}_2$-graduation in \eqref{so-triples}
with $\h=\sigma_{W,W}= \eespan \langle \sigma_{x,y}:x,y\in W\rangle$ (the maps $\sigma_{x,y}$ as in (\ref{sk-so})), and $\sigma_{v,W}= \eespan \langle \sigma_{v,x}:x\in W\rangle$, which is an irreducible module for $\h$ isomorphic to $W$. Thus from (\ref{condicionescasonosimpleII}), we have $\m=\h^\perp=\sigma_{v,W}$ which provides item (i) in the Lemma.

On the other hand, if $V$ is not irreducible as a module for $\h$, but
the restriction of $b$ to any irreducible $\h$-submodule is
trivial, by Weyl's theorem on complete reducibility, given
an irreducible submodule $W_1$ there is another irreducible
submodule $W_2$ with $b(W_1,W_2)\neq 0$. So, $W_1$ and $W_2$ are
isotropic, that is $b(W_i,W_i)= 0$, and contragredient modules, and $V=W_1\oplus W_2\oplus
(W_1\oplus W_2)^\perp$. Arguing as before, we may assume that $V=W_1\oplus W_2$.  Then our subalgebra $\h$ lies inside the subalgebra $\{f \in \so(V,b):
f(W_i)\subseteq W_i\}=\sigma_{W_1,W_2}$ and this contradicts the
maximality of $\h$, since $\sigma_{W_1,W_2}$ is contained properly in the
subalgebra $\sigma_{W_1,W_2}\oplus \sigma_{W_1,W_1}$.

\smallskip
Now, in case $V=V(\lambda)$ remains irreducible as a module for $\h$,
its dominant weight $\lambda$ relative to a Cartan subalgebra of $\h$ and a choice of a system $\Delta=\{\alpha_1, \dots, \alpha_n\}$ of simple roots, decomposes as $\lambda=\sum_{i=1}^n m_i\lambda_i$, where as in \cite{Hum72}, $\lambda_1,\ldots,\lambda_n$ denote the fundamental weights. Note that
$m_i=\langle \lambda,\alpha_i\rangle \geq 0$ is a non-negative integer for any $i$. Let $\alpha_i$ be a simple root which is not orthogonal to $\lambda$, that is $m_i\neq 0$. From
(\ref{copsk}), a copy of the irreducible module $V(2\lambda-\alpha_i)$ appears in
$\wedge^2V(\lambda)\cong \so(V)=\h \oplus \m$. In case $\h \cong
V(2\lambda-\alpha_i)$, we have that $2\lambda-\alpha_i$ is the highest root $\omega$ of $\h$, and hence $\omega+\alpha_i$ is twice a dominant weight $\lambda$ (while $\h$ being a proper subspace of $\Lambda^2(V(\lambda)$). A quick look at the Dynkin diagrams (see \cite{Hum72}) shows that the only possibilities are the ones that appear in items (ii-b) and (ii-c).

Otherwise we must assume that the highest root of $\h$ is not of the form $2\lambda-\alpha_i$ for some simple root $\alpha_i$ such that $\langle\lambda,\alpha_i\rangle\ne 0$. As $\so(V)$ has
exactly two irreducible components as a module for $\h$, there exists exactly one simple root $\alpha_i$ not orthogonal to $\lambda$. Hence
$\lambda=m_i\lambda_i$ with $m_i\ge1$ and
$\m=V(2m_i\lambda_i-\alpha_i)$ which provides item (ii-a).
\end{proof}

Following Lemma~\ref{estructura-so}, for any reductive and nonsymmetric decomposition $\so(V,b)=\h\oplus \m$ satisfying (a), (b) and (c) in (\ref{condicionescasonosimpleI}), the vector space $V$, considered as a module for $\h$  must be (nontrivial) irreducible with dominant weight of the form $m\lambda_i$,  $\lambda_i$ being a fundamental weight relative to some  system of simple roots $\Delta$ of $\h$. The irreducibility of $V$ allows us to endow this space with a structure of either a Lie triple system or an orthogonal triple system (see \cite[Section V]{Oku93}, \cite[Definition 4.1]{Eld06} or the Appendix subsection \ref{simplecticos} in this paper for the definition of the latter systems), in such a way that the subalgebra $\h$ becomes its inner derivation Lie algebra. In this way, the classification in the $\so$-case will follow from known results on these triple systems.

For an arbitrary reductive decomposition $\so(V,b)=\h\oplus \m$, by using the  isomorphism as modules for $\h$ in (\ref{sk-so}), we can define the map
\begin{equation}\label{proy-so}
\begin{array}{rclcl}
V\otimes V & \to  &\so(V, b)& \to& \h\\
x\otimes y & \mapsto &  \sigma_{x,y}&\mapsto & d_{x,y}
\end{array}
\end{equation}
where $d_{x,y}$ denotes the projection of the operator $\sigma_{x,y}$ onto $\h$, so
the subalgebra $\h$ can be written as $\h=\eespan \langle d_{x,y}:x,y\in V\rangle$.
Now, let us define the triple product  on $V$ given by
\begin{equation}\label{triple-so}
xyz:=d_{x,y}(z)
\end{equation} This  product satisfies the identities:
\begin{eqnarray}
& &\label{so1}xxz=0,\\
& &\label{so2} xy(uvw)=(xyu)vw+u(xyv)w+uv(xyw),\\
& &\label{so3} b(xyu,v)+b(u,xyv)=0,
\end{eqnarray}
for any $x,y,z,u,v\in V$.

Note that (\ref{so1}) is equivalent to the skew-symmetry of the operators $d_{x,y}$. Identity (\ref{so2}) tells us that the map given in (\ref{proy-so}) is a homomorphism of modules for $\h$ and
(\ref{so3}) follows from  $\h$ being a subalgebra of $\so(V,b)$. Moreover, since $d_{x,y}=xy\, \cdot$, we have that the subalgebra $\h$ is the inner derivation Lie algebra of the triple $V$:
\begin{equation}\label{indertriple}
\h=\eespan \langle xy \cdot:x,y\in V\rangle=\Inder V,
\end{equation}
and we get the following result:

\begin{lemma}\label{casotriple-ort-Lie}
Given a reductive decomposition $\so(V,b)=\h\oplus \m$ satisfying {\em (a), (b), (c)} in {\em (\ref{condicionescasonosimpleI})}, such that the vector space $V$ is irreducible as a module for $\h$, the vector space $V$ endowed with the triple product $xyz$ defined in {\em (\ref{triple-so})} is either a simple Lie triple system or a simple orthogonal triple system with associated bilinear form  $\xi b$ for some nonzero scalar $\xi$. Moreover, the subalgebra $\h$ satisfies the equation
\begin{equation}\label{inder-so}
\h=\eespan \langle xy\,\cdot\, :x,y\in V\rangle,
\end{equation}
and therefore coincides with the inner derivation Lie algebra of the corresponding triple system, and the subspace $\m$ is the orthogonal complement $\h^\perp$ to $\h$ relative to the Killing form of $\so(V, b)$.
\end{lemma}

\begin{proof}
First  we shall check that the vector space $\Hom_{\h}(\wedge^2V\otimes V, V)$ is $2$-dimensional. Since $V\cong V^*$ and $ \wedge^2V\cong \so(V)$,
\begin{equation}\label{casotriple-ort-Lie-iso}
\begin{array}{lcl}
\Hom_{\h}(\wedge^2V\otimes V, V)&\cong&\Hom_{\h}(\wedge^2V,V\otimes
V)\qquad\qquad\qquad\qquad\\
&\cong&\Hom_{\h}(\h,V\otimes V)\oplus\Hom_{\h}(\m,V\otimes V)
\end{array}
\end{equation}

Moreover, the irreducibility of $\h$ as a module for itself gives
\begin{equation}
\dimens \Hom_{\h}(\h,V\otimes V)=\dimens \Hom_{\h}(V\otimes V^*,\h)
\end{equation}
Lemma \ref{estructura-so} shows that  $V=V(m\lambda_i)$ (the irreducible module of dominant weight $m\lambda_i$) as a module for $\h$, so \cite[Theorem~1]{Fau85} proves that $\Hom_{\h}(\h,V\otimes V)$ is a one-dimensional vector space. From the different possibilities for  $\m$ described in (ii-a), (ii-b) and (ii-c) of Lemma \ref{estructura-so}, the same assertion holds for $\Hom_{\h}(\m,V\otimes V)$:
\begin{itemize}
\item[(ii-a)] $\m=V(m\lambda_i-\alpha_i)$.\null\quad\newline
    The assertion follows from (\ref{copsk}) and comments therein.
\item[(ii-b)] $\h \cong B_3$, $V=V(\lambda_3)$ and $\m=V(\lambda_1)$.\null\quad\newline
    The assertion follows since the tensor product decomposition
\begin{equation}\label{b3-so}
\begin{array}{rcl}
V(\lambda_3)\otimes V(\lambda_3)\cong V(2\lambda_3)\oplus V(\lambda_2)\oplus
V(\lambda_1)\oplus V(0)
\end{array}
\end{equation}
    contains only one copy of $\m$.
\item[(ii-c)] $\h \cong G_2$, $V=\m=V(\lambda_1)$.\null\quad\newline
    Again the tensor product decomposition
\begin{equation}\label{g2-so}
V(\lambda_1)\otimes V(\lambda_1)\cong V(2\lambda_1)\oplus
V(\lambda_2)\oplus V(\lambda_1)\oplus V(0)
\end{equation}
    contains only one copy of $\m$.
\end{itemize}

On the other hand, in a easy way we can get the following $\h$-module decomposition for the tensor product $\wedge^2V \otimes V$:

\begin{equation}\label{descomposicionbuena}
\wedge^2V\otimes V=\wedge^3 V \oplus S
\end{equation}
where $\wedge^3V$ embeds in $\wedge^2V\otimes V$ by means of $x\wedge y\wedge z\mapsto (x\wedge y)\otimes z+(y\wedge z)\otimes x+(z\wedge x)\otimes y$, and
$S=\eespan\langle x\wedge y\otimes z+z\wedge y\otimes x : x,y,z\in
V\rangle$. The nonzero $\h$-homomorphism $\varphi:S\to V$ given by:
\[
\begin{split}
\varphi((x\wedge y)\otimes z&+(z\wedge y)\otimes x)\\
 &=\sigma_{x,y}(z)+\sigma_{z,y}(x)=2b(x,z)y-b(y,z)x-b(y,x)z
\end{split}
\]
with $\sigma_{x,y}$ as in (\ref{sk-so}), provides the alternative decomposition
\begin{equation}\label{descomposicionbuenaII}
\wedge^2V\otimes V=\wedge^3V\oplus\mathrm{Ker}\varphi\oplus V
\end{equation}
and therefore
\begin{equation}\label{desc-hom}
\Hom_{\h}(\wedge^2V\otimes V, V)=\Hom_{\h}(\wedge^3V, V)\oplus\Hom_{\h}(\mathrm{Ker}\varphi, V)\oplus\Hom_{\h}(V, V)
\end{equation}

Since the dimension of the vector spaces $\Hom_{\h}(\wedge^2V\otimes V, V)$ and $\Hom_{\h}(V,V)$ is $2$ and $1$ respectively, either $\Hom_{\h}(\wedge^3V,V)=0$ or $\Hom_{\h}(S,V)= k\varphi$. In the first case, we have that the triple product $xyz$ defined in (\ref{triple-so}) must be trivial when restricted to $\wedge^3V$. Therefore this triple product satisfies the additional identity
\begin{equation}\label{lts2}
xyz+yzx+zxy=0
\end{equation}
for any $x,y,z\in V$.
Hence, from (\ref{so1}), (\ref{so2}) and (\ref{lts2}) we get that $(V,xyz)$ is a Lie triple system (see Subsection \ref{Lie} in the Appendix for the definition).  Moreover, as $\h=\Inder\ (V)$ and $V=V(m\lambda_i)$ is $\h$-irreducible, this triple system is simple.

Otherwise $\Hom_{\h}(S,V)= k\varphi$ holds, so the restriction of the triple product $xyz$ to $S$ give us the relationship
\begin{equation}\label{eq4.3}
xyz+zyx=\xi\varphi(x\wedge y\otimes z+z\wedge y\otimes x) =
2\xi b(x,z)y-\xi b(y,z)x-\xi b(y,x)z
\end{equation}
for some $\xi\in k$ and any $x,y,z\in V$. Moreover, let us show that $\xi$ must be nonzero. Assume opn the contrary that $\xi=0$, from (\ref{eq4.3}) we get
\begin{equation}\label{eq4.4}
xyz+zyx = 0
\end{equation}
for all $x,y,z\in V$ and the triple product is totally antisymmetric. Then the triple products
$\langle x_{\sigma}y_{-\sigma}z_{\sigma}\rangle=\sigma
x_{\sigma}y_{-\sigma}z_{\sigma}$ defined on the vector space pair $\mathcal{U}=(U^+,U^-)$ with $U^\sigma=V$ and $\sigma=\pm$ satisfy:
{\setlength\arraycolsep{2pt}
$$
\begin{array}{rcl}
\langle x_{\sigma}y_{-\sigma}z_{\sigma}\rangle &= &\sigma
x_{\sigma}y_{-\sigma}z_{\sigma}=-\sigma
z_{\sigma}y_{-\sigma}x_{\sigma}=-\langle
z_{\sigma}y_{-\sigma}x_{\sigma}\rangle
\end{array}
$$
and using (\ref{so1}) y (\ref{so2}),
$$
\begin{array}{rcl}
\langle x_{\sigma}y_{-\sigma}\langle
u_{\sigma}v_{-\sigma}w_{\sigma}\rangle\rangle & = & \sigma^2x_{\sigma}
y_{-\sigma}(u_{\sigma}v_{-\sigma}w_{\sigma})\\
& = &\sigma^2((x_{\sigma}y_{-\sigma}u_{\sigma})v_{-\sigma}w_{\sigma})-
u_{\sigma}(y_{-\sigma}x_{\sigma}v_{-\sigma})w_{\sigma}\\
& &+ u_{\sigma}v_{-\sigma}(x_{\sigma}y_{-\sigma}w_{\sigma}))\\
&=&\langle\langle
x_{\sigma}y_{-\sigma}u_{\sigma}\rangle v_{-\sigma}w_{\sigma}\rangle+
\langle u_{\sigma}\langle
y_{-\sigma}x_{\sigma}v_{-\sigma}\rangle w_{\sigma}\rangle\\
& & + \langle u_{\sigma}v_{-\sigma}\langle
x_{\sigma}y_{-\sigma}w_{\sigma}\rangle\rangle
\end{array}
$$
}Hence, $\mathcal{U}$ is an anti-Jordan pair for which the inner derivation operators are of the form
{\setlength\arraycolsep{1pt}
\begin{equation}\label{estasder}
\begin{array}{rl}
(D_+(x_+,y_-),D_-(y_-,x_+))&=(\langle x_+y_-\, \cdot\,\rangle, \langle y_-x_+\, \cdot\, \rangle)\\
&=(x_+y_-\, \cdot\,,-y_-x_+\, \cdot\,)
\end{array}
\end{equation}
}We note that the linearization $xyz=-yxz$ of (\ref{so1}) is equivalent to $D_{x,y}=-D_{y,x}$, thus $(D_+(x_+,y_-),D_-(y_-,x_+))=(D_{x_+,y_-},D_{x_+,y_-})$. This shows that the Lie algebra $\Inder\,  \mathcal{U}$ is isomorphic to $\h$ and therefore it is a simple Lie algebra. Moreover, since $V$ is an irreducible module for $\h$, $\mathcal{U}$ is a simple anti-Jordan pair (\cite[Proposition 1.2]{FauFe80}). But according to  Table \ref{simple-antijordan} in the Appendix, the inner derivation algebras of simple anti-Jordan pairs such that $\mathcal {U}^+=\mathcal {U}^-$ are not simple. Hence $\xi\neq 0$.

Now equation (\ref{eq4.3}) with $z=y$ and ( \ref{so1}) give the identity
\begin{equation}\label{so4}
xyy=\xi b(x,y)y-\xi b(y,y)x
\end{equation}
for any $x,y\in V$,
which together with (\ref{so1}), (\ref{so2}) and (\ref{so3}) (see Subsection \ref{simplecticos} in the Appendix) prove that the vector space $V$ is an orthogonal triple system under the triple product $xyz$  and the symmetric bilinear form $\xi b(x,y)$. Since the form is nondegenerate, $V$ is a simple orthogonal triple system (\cite[Proposition 4.4]{Eld06}) and the subalgebra $\h$ is its inner derivation Lie algebra. The last assertion ($\m=\h^\perp$) follows from condition (d) in
(\ref{condicionescasonosimpleII})
\end{proof}

Now we can formulate the main result for the generic $\so$-case:

\begin{theorem}\label{descripcion-so}
Let $(\m, a\cdot b, [a,b,c])$ be an irreducible LY-algebra of generic type and standard enveloping Lie algebra of type $\so$. Then either:
\begin{itemize}
\item[(i)] There is a vector space $V$ of dimension $\geq 5$, endowed with a nondegenerate symmetric bilinear form $b$ such that $\m$ is, up to isomorphism, the simple Lie triple system defined on $V$ with triple product $[u,v,w]=b(u,w)v-b(v,w)u$. In particular, the binary product $u\cdot v$ is trivial.

\item[(ii)] Up to isomorphism, $\m$ coincides with the space $\cO_0$ of zero trace octonions with binary and ternary products $a\cdot b=ab-ba=[a,b]$ and
    $[a,b,c]=2\bigl([[a,b],c]-3\bigl((ac)b-a(cb)\bigr)\bigr)$ for any $a,b,c\in\cO_0$, where $ab$ denotes the multiplication in $\cO$.

\item[(iii)] There is a simple Lie triple system $T$ endowed with a nondegenerate symmetric bilinear form $b$ of one of the following types:
    \begin{enumerate}
    \item[(iii.a)] a simple Lie algebra of type different from $A$ with its natural triple product $[xyz]=[[x,y],z]$ endowed with its Killing form,
    \item[(iii.b)] the subspace of zero trace elements of a simple Jordan algebra of degree $\geq 3$, not isomorphic neither to $\Mat_n(k)^+$ ($n\geq 3$) nor to $\cH_4(k)$, with its triple product $[xyz]=(x,z,y)=(x\circ z)\circ y-x\circ(y\circ z)$ (where $x\circ y$ denotes the multiplication in the Jordan algebra), endowed with the nondegenerate bilinear form given by its generic trace,
    \item[(iii.c)] the Lie triple systems attached to the exceptional symmetric pairs $(F_4,B_4)$, $(E_6,C_4)$, $(E_7,A_7)$ or $(E_8,D_8)$, endowed with the nondegenerate bilinear form given by the restriction of the Killing form of the ambient Lie algebra,
    \end{enumerate}
    such that, up to isomorphism, $\g(\m)=\so(T,b)$, $\h=D(\m,\m)=\Der T=[TT.]=\eespan\langle [xy.]:x,y\in T\rangle$. Here the LY-algebra $\m$ appears as the orthogonal complement to $\h$ in $\g(\m)$ relative to the Killing form, with the binary and ternary products in \eqref{eq:binter}.
\end{itemize}

There are no isomorphisms among the LY-algebras in the different items above.
\end{theorem}

\begin{proof} From Lemmas \ref{estructura-so} and \ref{casotriple-ort-Lie} we know that  either
$\m=\sigma_{v,W}=\eespan\langle \sigma_{v,x}: x \in W\rangle$, where $\sigma_{v,x}=b(v,\cdot)x-b(x,\cdot)v$, inside the symmetric decomposition $\so(kv\oplus W, b)=\sigma_{W,W}\oplus \sigma_{v,W}$ where $b$ is a nondegenerate symmetric form with $b(v,W)=0$ or $\m=\h^\perp$ where $\h$ is the linear Lie algebra of inner derivations related to either a simple Lie triple system or a simple orthogonal triple system which is irreducible as a module for its inner derivation algebra.

In the first case, since $[\sigma_{x,y},\sigma_{a,b}]=\sigma_{\sigma_{x,y}(a),b}+\sigma_{a,\sigma_{x,y}(b)}$ for any $x,y,a,b\in V$, it follows that
\begin{equation}
[[\sigma_{v,x}, \sigma_{v,y}],\sigma_{v,z}]=b(v,v)\sigma_{v, \sigma_{x,y}(z)},
\end{equation}
for any $x,y,z\in W$. Besides, since the ground field $k$ is assumed to be algebraically closed, one may take $v$ with $b(v,v)=1$. Hence $\m=\sigma_{v,W}$ can be identified to $W$ with trivial binary product and triple product given by $[xyz]=\sigma_{x,y}(z)$, thus obtaining the situation in item (i)

Otherwise,  we must look for either simple Lie triple systems or simple orthogonal  triple systems $(V,xyz)$ such that:
\begin{itemize}
\item[(1)] the inner derivation algebra $\Inder V$ is simple,
\item[(2)] $V=V(m\lambda_i)$ is an irreducible module for $\Inder\,V$ with dominant weight $m$-times a fundamental weight $\lambda_i$, and $\dimens V \geq 5$,
\item[(3)] $\wedge^2V$ decomposes as a sum of two irreducible modules.
\end{itemize}

Since isomorphic irreducible orthogonal or Lie triple systems provide isomorphic LY-algebras, we need to check the previous conditions in the classifications, up to isomorphisms, of such  systems  given in \cite[Theorem 4.7]{Eld06}, \cite[Table I]{Fau80} and \cite[Table III]{Fau85}, which  are outlined in the Appendix of this paper: Tables \ref{ltsadjuntos-so-sl-excep} and \ref{simple-ortogonales}. Then, using  Table \ref{simple-ortogonales} in the Appendix, and restrictions (1)-(2)-(3), we get the following possibilities for the triple $(V,\Inder V, V(m\lambda_i))$  for orthogonal triple systems:  the $G$-type triple system defined on the $7$-dimensional space $\cO_0$  of zero trace octonions with triple $(\cO_0, G_2, V(\lambda_1))$ and the $F$-type triple systems defined on a $8$-dimensional vector space $V$ having a $3$-fold vector cross product with ternary description $(V, \so_7(k),V(\lambda_3))$. The respective decompositions of  $\wedge^2 V$ as a module for $\h=\Inder\, V$ are the following:

\begin{description}
\settowidth{\labelwidth}{XX}%
\setlength{\leftmargin}{30pt}

\item[$(\cO_0, G_2, V(\lambda_1))$]\null\quad\newline
    $\wedge^2V(\lambda_1)=V(\lambda_1)\oplus V(\lambda_2)$ as a module for $\Inder\, \cO_0=G_2$.

\item[$(V, \so_7(k), V(\lambda_3))$]\null\quad\newline
    $\wedge^2V(\lambda_3)=V(\lambda_1)\oplus V(\lambda_2)$ as a module for $\Inder\, V\simeq\so_7(k)$.
\end{description}

Following \cite{Eld06}, the orthogonal triple system of $G$-type satisfies that $\h=\Inder\,V$ is the simple Lie algebra of type $G_2$ given by the Lie algebra of derivations of $\cO$, considered as a subalgebra of $\so(\cO_0,n)$, where $n$ denotes the norm of the octonion algebra. Then $\m$ is the orthogonal complement of $\h$ relative to the Killing form of $\g=\so(\cO_0,n)$. But $\so(\cO_0,n)$ decomposes as
\[
\so(\cO_0,n)=\Der(\cO)\oplus\ad_{\cO_0},
\]
(see \cite[Chapter III, \S 8]{Sch} or \cite{EldMyung93}), where $\ad_x(y)=[x,y]=xy-yx$. Since $\ad_{\cO_0}$ is irreducible as a module for $\Der(\cO)$ it turns out that $\ad_{\cO_0}$ is necessarily the orthogonal complement to $\h=\Der(\cO)$ relative to the Killing form of $\so(\cO_0,n)$. Besides for any $x,y\in\cO_0$ we have:
\[
\begin{split}
[\ad_x,\ad_y]&=[L_x-R_x,L_y-R_y]\\
 &=[L_x,L_y]+[R_x,R_y]-[R_x,L_y]-[L_x,R_y]\\
 &=D_{x,y}-3[L_x,R_y],
\end{split}
\]
where $L_x$ and $R_x$ denote the left and right multiplications by $x$ in $\cO$, and where $D_{x,y}=[L_x,L_y]+[L_x,R_y]+[R_x,R_y]=\ad_{[x,y]}-3[L_x,R_y]$ is the inner derivation of $\cO$ generated by the elements $x$ and $y$. Therefore
\[
[\ad_x,\ad_y]=-\ad_{[x,y]}+2D_{x,y},
\]
and hence, if we identify $\m=\ad_{\cO_0}$ with $\cO_0$ by means of $x\mapsto -\ad_x$, the binary and ternary multiplications in \eqref{eq:binter} are given by:
\[
\begin{split}
x\cdot y&=[x,y],\\
[xyz]&=2D_{x,y}(z),
\end{split}
\]
for any $x,y,z\in\cO_0$ and we obtain item (ii).

For $F$-type orthogonal triple systems (see \cite{Eld06} and \cite{EldKaO}), we have that $V=\cO$ with bilinear form $b(x,y)=\alpha n(x,y)$ where $n(x,y)$ is as in the previous paragraph, the triple product is given by $xyz=(x{\bar y})z+4b(x,z)y-4b(y,z)x-b(x,y)z$ and $\Inder V=\eespan \langle xy\, \cdot\rangle$. In this case $\Inder V$ is a Lie algebra of type $B_3$ that can be described as $\Inder V=\eespan\langle D_{x,y}, L_x+2R_x: x,y \in \cO_0\rangle$. Moreover, the automorphism $\theta$ of $\so(\cO,n)\cong \so_8(k)$ given by $L_x \mapsto L_x+R_x$ and $R_x \mapsto -R_x$ (see \cite[Chapter III,\S 8]{Sch} or \cite[Theorem 3.2]{Eld00}) makes $\theta(\Inder V)=\eespan\langle D_{x,y}, L_x-R_x=\ad_x: x,y \in \cO_0\rangle=\{f \in \so(V): f(1)=0\}=\so(\cO_0, n)$, a Lie algebra for which $\cO_0$ is irreducible and orthogonal to $k 1$ which is a one-dimensional and $\theta(\Inder V)$-invariant subspace. Hence, the LY-algebra $\m=(\Inder V)^\perp$ is isomorphic to $\sigma_{1,\cO_0}=\eespan \langle \sigma_{1,x}=n(1,\cdot)x-n(x,\cdot)1: x\in \cO_0\rangle$ obtained  as in item (i). Hence nothing new appears here.

For simple Lie triple systems, \cite[Table I]{Fau80} presents the complete classification of such triple systems encoded through affine Dynkin diagrams. Using this classification, in \cite[Table III]{Fau85} a complete list of all simple Lie triple systems which are irreducible for $\Inder V$ is given. Table III also provides the inner derivation algebra $\Inder V$, and the  structure for each irreducible Lie triple $V$ as a module for $\Inder(V)$. The results in \cite{Fau80} and \cite{Fau85} are displayed on Table \ref{ltsadjuntos-so-sl-excep} in the Appendix. Table \ref{ltsadjuntos-so-sl-excep} provides the dominant weights for the different irreducible Lie triple systems $V$ as well as the $\wedge^2V$-decomposition of those triples with simple inner derivation Lie algebra. So using Table \ref{ltsadjuntos-so-sl-excep} under the restrictions (1)-(2)-(3), we arrive at the possibilities described below. We also note that in all the cases, $V$ is a contragredient and irreducible module for $\Inder(V)$, and hence there exists a unique $\Inder V$-invariant  form $b$ on $V$ up to scalars, which is either symmetric or skew-symmetric. The $\wedge^2V$-decomposition in Table \ref{ltsadjuntos-so-sl-excep} proves that the form $b$ is always symmetric for irreducible triple systems with simple inner derivation algebra.

\begin{description}
\settowidth{\labelwidth}{XX}%
\setlength{\leftmargin}{30pt}

\item[$(\cL\times \cL, \cL)$ ]\null\quad\newline
    Here $\cL$ represents a Lie algebra of type different from $A_n, n\geq 1$, but considered as a triple system by means of the product $[xyz]=[[x,y],z]$, and they are endowed with the Killing form. These triple systems are the so called adjoint in \cite{BEM09}.

\item[$(A_{2n}, B_n)_{n\geq 1}$, $(A_{2n-1}, C_n)_{n\geq 3}$, $(A_{2n-1}, D_n)_{n\geq 3}$]\null\quad\newline
    Following Theorem \ref{descripcion-sl}, these triple systems correspond to the space of zero trace elements in a simple Jordan algebra $\cJ$ of type $B$ or $C$ and degree $n \geq 3$ with triple product given by the associator. Then, since $(y,z,x)=[L_x,L_y](z)$ ($L_x$ denotes the multiplication by $x$ in the Jordan algebra), we have that $\Inder \cJ_0=\eespan \langle [L_x,L_y]:x,y \in \cJ_0  \rangle=\Der \cJ$ and $b$ is the generic trace.

\item[$(F_4, B_4)$, $(E_{6}, F_4)$, $(E_6, C_4)$, $(E_7, A_7)$, $(E_8,D_8)$]\null\quad\newline
    Table \ref{ltsadjuntos-so-sl-excep} shows that all the simple Lie triple systems with exceptional simple standard enveloping Lie algebra and simple inner derivation Lie algebra work. The triple system related to the pair $(E_{6}, F_4)$ consists of the zero trace elements of the exceptional simple Jordan algebra (Albert algebra) with associator as triple product (see \cite{JacJA}).
\end{description}

On the other hand, Table \ref{ltsadjuntos-so-sl-excep} shows that the only symmetric decompositions with standard enveloping Lie algebra of type $\so$ are given by symmetric pairs of type $(\so_n(k),\so_{n-1}(k))$, up to isomorphisms. It is easy to check that none of the reductive pairs related with the LY-algebras described in items (ii) and (iii)  are of this form. So the binary and ternary products of the corresponding LY-algebras are not trivial.

Finally, the restriction on item (iii.b) on the Jordan algebra not being isomorphic to $\cH_4(k)$ is due to the fact that for this Jordan algebra, the associated subalgebra $\h$ is $\so_4(k)$ which is not simple.
\end{proof}

\begin{remark}
The case in item (iii.b) of the previous Theorem corresponding to the Jordan algebra $\cH_3(k)$ satisfies that its enveloping algebra is $\so_5(k)$, which is isomorphic to $\frsp_4(k)$. Hence this case will appear too in the next section (Item (i) of Theorem \ref{teo-sp}). Therefore, we may assume that the Jordan algebra in item (iii.b) above is not isomorphic to $\cH_3(k)$. This will be done in our final Table \ref{clasificacion generico}.
\end{remark}


\section{Symplectic case}
For LY-algebras of Generic Type and standard enveloping Lie algebra a (simple) symplectic Lie algebra $\frsp(V,b)$, we will follow a similar procedure to that used in the special and orthogonal cases. In the symplectic case, $V$ is an even-dimensional vector space endowed with a nondegenerate skew-symmetric form $b$.

Given a suitable reductive decomposition $\frsp(V)=\h\oplus \m$, we may view $V$ and $\frsp(V)$ as modules for $\h$.  The map $x\cdot y \mapsto \gamma_{x,y}=b(x,\cdot)y+b(y,\cdot)x$ provides an isomorphism of $\h$-modules:
\begin{equation}\label{sym-sp}
S^2(V)\cong \frsp(V)
\end{equation}
where $S^2(V)$ is the second symmetric power of $V$. This isomorphism and the following easy Lemma on representations of Lie algebras are used in an essential way along this section:

\begin{lemma}\label{lema2-so} Let $\mu_1$ and $\mu_2$ be two dominant weights of a simple Lie algebra (relative to a fixed system of simple roots). Then the modules  $\wedge^2(V(\mu_1))$ and $S^2(V(\mu_2))$ are isomorphic if and only if one of the following holds:
\begin{enumerate}
\item[(i)] $\h$ is a simple Lie algebra of type $A_1$, $\mu_1=k\lambda_1$ and $\mu_2=(k-1)\lambda_1$ with $k\ge 1$,
\item[(ii)] $\h$ is a simple Lie algebra of type $B_2$, $\mu_1=\lambda_1$ and $\mu_2=\lambda_2$.
\end{enumerate}
\end{lemma}
\begin{proof}

For a given simple root $\alpha_i$ non-orthogonal to $\mu_1$, the weight $2\mu_1-\alpha_i$ is maximal in the set of weights of the module $\wedge^2V(\mu_1)$ relative to the usual partial order where $\lambda>\mu$ if $\lambda -\mu$ is a sum of positive roots. Since $2\mu_2$ is the only maximal weight for $S^2V(\mu_2)$, we have $2\mu_1-\alpha_i=2\mu_2$ and there is a unique simple root $\alpha_i$ not orthogonal to $\mu_1$. The last assertion implies $\mu_1=k\lambda_i$, $k\ge 1$ and $\alpha_i=2(\mu_1-\mu_2)$. Now it is easy to check that the only possibilities are the following (see \cite[Chapter III, Section 11]{Hum72}):
\begin{itemize}
\item $A_1$ with $\alpha_i=\alpha_1=2\lambda_1$, which implies item (i) in Lemma.

\item $B_2$ and $\alpha_i=\alpha_1=2(\lambda_1-\lambda_2)$:\null\quad\newline
    In this case, $\mu_1=k\lambda_1$ and $\mu_2=(k-1)\lambda_1+\lambda_2$. But computing the dimension of the corresponding irreducible modules $V(k\lambda_1)$ and $V((k-1)\lambda_1+\lambda_2)$, we get $k=1$ as the only possibility, thus item (ii) in Lemma follows.
\item $\h=C_n$, $n\ge 3$ and $\alpha_i=\alpha_n=2(-\lambda_{n-1}+\lambda_n)$:\null\quad\newline
    Then    $\mu_1=k\lambda_n$ and $\mu_2=\lambda_{n-1}+(k-1)\lambda_n$. But the  formula
$$\frac{\dim V(k\lambda_n)}{\dim
V((k-1)\lambda_n+\lambda_{n-1})}=\frac{2k+n+1}{2kn}
$$
    implies that $\dim V(k\lambda_n)< \dim V((k-1)\lambda_n+\lambda_{n-1})$ except for $n=3$ and $k=1$. So, the only possibility for both modules to be isomorphic is $\mu_1=\lambda_3$ and $\mu_2=\lambda_2$. But then $\dimens\wedge^2V(\lambda_3)=91<\dimens S^2V(\lambda_2)=105$ and therefore this situation does not hold.
\end{itemize}
The converse is easily checked by using the Clebsch-Gordan formula and the isomorphism between the $B_2$-type Lie algebra $\so_5(k)$ and the $C_2$-type $\frsp_4(k)$.
\end{proof}

Recall that given an irreducible module $V(\lambda)$ for a dominant weight $\lambda$ of a simple Lie algebra, the dual module $V(\lambda)^*$ is isomorphic to $V(-\sigma\lambda)$, where $\sigma$ is the element of the Weyl group sending the given system of simple roots to its opposite (see \cite[\S 21, Exercise 6]{Hum72}). We will write $-\sigma\lambda=\lambda^*$ and will say that the dominant weight $\lambda$ is self-dual in case $\lambda=\lambda^*$, that is, in case $V(\lambda)$ is a self-dual module.

\begin{lemma}\label{estructura-sp}
Let $\frsp(V)=\h\oplus\m$ be a reductive decomposition satisfying {\em
(a)}, {\em (b)} and {\em (c)} in {\em
(\ref{condicionescasonosimpleI})}.  Then  $\dimens V \geq 4$  and as a module for  $\h$, $V=V(\lambda)$ is irreducible and either its dominant weight $\lambda$ is a fundamental and self-dual weight or $\h=A_1$, $\dimens V=4$ and $\lambda=3\lambda_1$. In any case, $\m$ is an irreducible module for $\h$ whose dominant weight is $2\lambda$.
\end{lemma}

\begin{proof}
In case $V$ is reducible as a module for $\h$, the arguments in the proof of Lemma \ref{estructura-so} show that the vector space $V$ can be decomposed as an orthogonal sum, $V=W\oplus W^\perp$ with $W$ irreducible and nontrivial and $\h=\frsp(W)\oplus \frsp (W^\perp)$ which is not a simple Lie algebra. Hence $V$ must be irreducible and the assertion on its dimension  is clear.  In the sequel let $\lambda$ be the dominant weight of the irreducible self-dual module $V$.

On the other hand,  (\ref{copsim}) and (\ref{sym-sp}) show that $V(2\lambda)$ appears as a submodule in $S^2 V\cong\frsp(V)= \h\oplus\m$, thus as modules over $\h$ either $\h$ or $\m$ is isomorphic to $V(2\lambda)$.  In the first case the only possibility is that $2\lambda=2\lambda_1$ for the simple Lie algebra of type $C_n$. This implies $\frsp(V)=\frsp(V(\lambda_1))=\h$, which is not possible. Hence $\m$ is irreducible with $2\lambda$ as dominant weight.

Now let $\lambda=\sum m_i \lambda_i$ the decomposition of $\lambda$ as sum of fundamental weights $\lambda_i$, and assume that $\lambda$ is not fundamental. Then $\lambda$ can be decomposed in two different ways:
\begin{equation}\label{lambdadesc}
\begin{array}{ll}
\textrm{a)}&\lambda=\lambda_i+\lambda_i^* \\
&\mathrm {with}\ \lambda_i \  \textrm {fundamental\ and\ non self-dual}\\

\textrm{b)}&\lambda=\lambda'+\lambda'' \\
&\mathrm {with}\ \lambda', \lambda''  \ \textrm{nonzero\ and\ self-dual\ dominant\  weights}
\end{array}
\end{equation}

Suppose $\lambda=\lambda_i+\lambda_i^*$ and note that $\prod(\lambda_i^*)=\{-\mu:\mu\in\prod(\lambda_i)\}$ is the set of weights for the module $V(\lambda_i)^*=V(\lambda_i^*)$. (As in \cite{Hum72} $\prod(\lambda)$ denotes the set of weights of $V(\lambda)$.) Since $f_{\mu}(v_{-\mu})\ne 0$ in case $v_{-\mu}\in V(\lambda_i)$ and $f_{\mu}\in V(\lambda_i^*)$ are $(\pm \mu)$-weight vectors, the symmetric $\h$-invariant form
\begin{equation}\label{b1-sp}
\begin{array}{rcl}
b:(V(\lambda_i)\otimes V(\lambda_i)^*)\otimes(V(\lambda_i)\otimes V(\lambda_i)^*) &
\to & F \\
(v_1\otimes f_1)\otimes (v_2\otimes f_2) & \mapsto & f_1(v_2)f_2(v_1)
\end{array}
\end{equation}
satisfies $b(v_{\lambda_i}\otimes f_{\lambda_i^*},v_{-\lambda_i^*}\otimes f_{-\lambda_i})\ne 0$. As the copy of $V(\lambda)$ generated by $v_{\lambda_i}\otimes f_{\lambda_i^*}$ that appears in $V(\lambda_i)\otimes V(\lambda_i)^*$ contains too the element $v_{-\lambda_i^*}\otimes f_{-\lambda_i}$, the symmetric form $b$ induces a symmetric and nonzero $\h$-invariant form on $V(\lambda)$, but this is not possible: because of the irreducibility of $V$, up to scalars there is exactly one $\h$-invariant form on $V$, which must be skew-symmetric.

Hence from (\ref{lambdadesc}) $\lambda$ decomposes as $\lambda=\lambda'+\lambda''$. Then,  the self-dual modules $V(\lambda')$ and $V(\lambda'')$ are endowed with nondegenerate and $\h$-invariant forms $b_1(x',y')$ and $b_2(x'',y'')$. From $b_1$ and $b_2$, we can define on the tensor product $V(\lambda')\otimes V(\lambda'')$ the $\h$-invariant form
\begin{equation}\label{b2-sp}
\begin{array}{rcl}
\hat b:V(\lambda')\otimes V(\lambda'')\otimes V(\lambda')\otimes V(\lambda'') &
\to & F \\
v'_1\otimes v''_1\otimes v'_2\otimes v''_2 & \mapsto &
b_1(v'_1,v'_2)b_2(v''_1,v''_2)
\end{array}
\end{equation}
which satisfies $\hat b(v_{\lambda'}\otimes v_{\lambda''}, v_{-\lambda'}\otimes v_{-\lambda''})\ne 0$. Now, a copy of $V(\lambda)^{\otimes^2}$ appears in $(V(\lambda')\otimes V(\lambda''))^{\otimes^2}$ so, $\hat b$ defines a nonzero and $\h$-invariant form on $V(\lambda)^{\otimes^2}$ which must be skew-symmetric. Consequently and without loss of generality, we can assume that $b_1$ is symmetric and $b_2$ is skew-symmetric. Now let $c'$ and $c''$ be the module homomorphisms given by

\begin{equation}
\begin{array}{rcl}
c':V(\lambda')\otimes V(\lambda'')\otimes V(\lambda')\otimes V(\lambda'') &
\to & V(\lambda'')\otimes V(\lambda'') \\
v'_1\otimes v''_1\otimes v'_2\otimes v''_2 & \mapsto &
b_1(v'_1,v'_2)v''_1\otimes v''_2
\end{array}
\end{equation}
and
\begin{equation}
\begin{array}{rcl}
c'':V(\lambda')\otimes V(\lambda'')\otimes V(\lambda')\otimes V(\lambda'') &
\to & V(\lambda')\otimes V(\lambda') \\
v'_1\otimes v''_1\otimes v'_2\otimes v''_2 & \mapsto &
b_2(v''_1,v''_2)v'_1\otimes v'_2
\end{array}
\end{equation}
Since $b_1$ is symmetric, we have
{\setlength\arraycolsep{2pt}
$$
\begin{array}{rl}
c'(v_{\lambda'}\otimes v_{\lambda''}\otimes v_{-\lambda'}\otimes v_{-\lambda''} & +
v_{-\lambda'}\otimes v_{-\lambda''}\otimes v_{\lambda'}\otimes v_{\lambda''})\\
&=b_1(v_{\lambda'}, v_{-\lambda'})(v_{\lambda''}\otimes v_{-\lambda''} +
v_{-\lambda''}\otimes v_{\lambda''})\ne 0
\end{array}
$$
}Hence, as the symmetric modules $S^2(V(\lambda))$ and $S^2(V(\lambda''))$ are generated by the vectors $v_{\lambda'}\otimes v_{\lambda''}\otimes v_{-\lambda'}\otimes v_{-\lambda''} +
v_{-\lambda'}\otimes v_{-\lambda''}\otimes v_{\lambda'}\otimes v_{\lambda''}$ and $v_{\lambda''}\otimes v_{-\lambda''} +v_{-\lambda''}\otimes v_{\lambda''}$ respectively, we get  a copy of the second symmetric power $S^2(V(\lambda''))$ inside $S^2(V(\lambda))$. Moreover, since $\dimens V(\lambda) >\dimens V(\lambda'')$, we have $S^2(V(\lambda''))\neq S^2(V(\lambda))$. In this way using (\ref{sym-sp}), $S^2(V(\lambda''))$ appears as a proper submodule inside $\frsp(V)=\h\oplus\m$. Similar arguments for  $c''$ and the skew symmetric $b_2$ yields
{\setlength\arraycolsep{2pt}
$$
\begin{array}{rl}
c''(v_{\lambda'}\otimes v_{\lambda''}\otimes v_{-\lambda'}\otimes v_{-\lambda''} & +
v_{-\lambda'}\otimes v_{-\lambda''}\otimes v_{\lambda'}\otimes v_{\lambda''})\\
&=b_2(v_{\lambda''}, v_{-\lambda''})(v_{\lambda'}\otimes v_{-\lambda'} -
v_{-\lambda'}\otimes v_{\lambda'})\ne 0,
\end{array}
$$
}and therefore the second alternating power $\wedge^2(V(\lambda'))$ also appears properly on the module decomposition of $S^2(V(\lambda))\cong \h\oplus\m$. Since $\h$ is contained in both $\so(V(\lambda'),b_1)\simeq\wedge^2(V(\lambda'))$ and in $\frsp(V(\lambda''),b_2)\simeq S^2(V(\lambda''))$, and $\m$ is irreducible, we get
\begin{equation}
\so(V(\lambda'))\cong\wedge^2(V(\lambda'))\cong\h\cong
S^2(V(\lambda''))\cong\frsp(V(\lambda''))
\end{equation}
Then Lemma \ref{lema2-so} shows that  either $\h$ is a simple Lie algebra of type $A_1$,
$\lambda'=2\lambda_1$ and $\lambda''=\lambda_1$; or $\h$ is simple of type $B_2$ and therefore $\lambda'=\lambda_1$ and $\lambda''=\lambda_2$. The latter possibility does not work from the dimensionality of the different modules involved: $\dim V=\dim V(\lambda_1+\lambda_2)=16$, $\dim \h= \dimens B_2=10$ and $\dim \m=\dim V(2\lambda_1+2\lambda_2)=81$, but $\dim \frsp(V)=136>\dim \h+\dim \m=91$. Hence, in case $\lambda$ is not fundamental, $\h$ is of type $A_1$ with $\lambda=\lambda'+\lambda''=3\lambda_1$ and this completes the proof.
\end{proof}
\bigskip

Lemma \ref{estructura-sp} shows that the irreducible LY-algebras which appear inside reductive decompositions $\frsp(V)=\h\oplus \m$ satisfying (a), (b) and (c) in (\ref{condicionescasonosimpleI}) are given by simple and maximal linear subalgebras $\h$ with natural action on $V$ given by a fundamental and self-dual dominant weight except for $\frsp_4(k)\cong \so_5(k)$. This  allows us to endow $V$ with a structure of either a symplectic triple system(see \cite{YamAs} and \cite[Definition 2.1]{Eld06} for a definition) or an anti-Lie triple system (see \cite{FauFe80}), such that  $\h$ becomes its inner derivation algebra. In this way, the classification in the $\frsp$-case will follow from known results on these triple systems.

For an arbitrary reductive decomposition $\frsp(V,p)=\h\oplus \m$, the $\h$-module isomorphism in (\ref{sym-sp}) allows us to define the map
\begin{equation}\label{sp-map}
\begin{array}{rclcl}
V\otimes V & \to  &\frsp(V,b)& \to& \h\\
x\otimes y & \mapsto &  \gamma_{x,y}&\mapsto & d_{x,y}
\end{array}
\end{equation}
where $d_{x,y}$ denotes the projection of $\gamma_{x,y}=b(x,\cdot)y+b(y,\cdot)x$ onto $\h$, so the subalgebra $\h$ appears as $\h=\eespan\langle d_{x,y}: x,y \in V\rangle$. Using these projections $d_{x,y}$, we define the triple product on $V$
\begin{equation}\label{triple-sp}
xyz:=d_{x,y}(z)
\end{equation}
which satisfies the following identities:
\begin{eqnarray}
& &\label{sp1}xyz=yxz\\
& &\label{sp2} xy(uvw)=(xyu)vw+u(xyv)w+uv(xyw)\\
& &\label{sp3} b(xyu,v)+b(u,xyv)=0
\end{eqnarray}
for $x,y,z \in V$. Identity (\ref{sp1}) is equivalent to the symmetry of the operators $d_{x,y}$. Identity (\ref{sp2}) states that (\ref{sp-map}) is an $\h$-module homomorphism and
(\ref{sp3}) follows because $\h$ is a subalgebra of $\frsp(V,b)$. Moreover, since $d_{x,y}=xy \, \cdot$, the subalgebra $\h$ becomes the inner derivation algebra of the triple $(V,xyz)$, so that
\begin{equation}\label{indertriplesp}
\h=\eespan \langle xy\, \cdot:x,y\in V\rangle=\Inder\  (V).
\end{equation}
Then, we have the following result, which is parallel to Lemma \ref{casotriple-ort-Lie}.

\begin{lemma}\label{casotriple-sim-aLie}
Given a reductive decomposition $\frsp(V,b)=\h\oplus \m$   satisfying {\em (a), (b), (c)} in {\em (\ref{condicionescasonosimpleI})}, the vector space $V$ endowed with the triple product $xyz$ defined in {\em (\ref{triple-sp})} is either a simple anti-Lie triple system of classical type or a simple symplectic triple system with associated symmetric form $\xi b$ for some nonzero scalar $\xi$. Moreover, the subalgebra $\h$ satisfies the equation
\begin{equation}\label{inder-sp}
\h=\eespan \langle xy\,\cdot\, :x,y\in V\rangle
\end{equation}
and therefore coincides with the inner derivation algebra of the corresponding triple system, and the subspace $\m$ is the orthogonal complement $\h^\perp$ to $\h$ relative to the Killing form of $\frsp(V,b)$.
\end{lemma}

\begin{proof} Since the triple product (\ref{triple-sp}) belongs to $\Hom_{\h}(S^2V\otimes V, V)$, we will describe the previous vector space in order to get the different possible products. Following Lemma \ref{estructura-sp}, $V$ is an irreducible and self-dual module for $\h$. Also (\ref{sym-sp}) gives $S^2V\cong\frsp(V)=\h\oplus\m$, so
{\setlength\arraycolsep{1pt}
\begin{equation}\label{prop2-sp-iso}
\begin{array}{rcl}
\Hom_{\h}(S^2V\otimes V, V)& \cong & \Hom_{\h}(S^2V,V\otimes V)\\
& \cong & \Hom_{\h}(\h,V\otimes V)\oplus\Hom_{\h}(\m,V\otimes V)
\end{array}
\end{equation}
}Then, using \cite[Theorem 1]{Fau80} and the dimension equality
$$
\dimens \Hom_{\h}(\h,V\otimes V)=\dimens \Hom_{\h}(V\otimes V,\h)
$$
the first summand in (\ref{prop2-sp-iso}) is one-dimensional. The same is deduced from the Clebsch-Gordan formula for the second summand in case $\h$ is of type $A_1$ and $V=V(3\lambda_1)$, as then we have $V(3\lambda_1)\otimes V(3\lambda_1)\cong V(6\lambda_1)\oplus V(4\lambda_1)\oplus V(2\lambda_1)\oplus V(0)$. Otherwise  following Lemma \ref{estructura-sp}, $V=V(\lambda_i)$ with $\lambda_i$ fundamental and $\m=V(2\lambda_i)$,  so the result follows from (\ref{copsim}). Hence, $\Hom_{\h}(S^2V\otimes V, V)$ is always a  two dimensional vector space.

On the other hand, $S^2V\otimes V$ can be decomposed as the module sum:
\begin{equation}\label{descs3-sp-S}
S^2V\otimes V=S^3V \oplus \eespan \langle (x\otimes y+y\otimes x)\otimes z-(z\otimes y+y\otimes z)\otimes x \rangle
\end{equation}
Write $S=\eespan\langle (x\otimes y+y\otimes x)\otimes z-(z\otimes y+y\otimes z)\otimes x: x,y,z, \in V \rangle$, then we can consider the nonzero $\h$-module homomorphism $\varphi:S\to V$ given by
{\setlength\arraycolsep{2pt}
$$
\begin{array}{l}
\varphi((x\otimes y + y\otimes x)\otimes z-(z\otimes
y+y\otimes z)\otimes x)\\ \qquad \qquad\qquad \qquad= \gamma_{x,y}(z)-\gamma_{z,y}(x)
=2b(x,z)y+b(y,z)x-b(y,x)z
\end{array}
$$
}where $\gamma_{x,y}$ is given in (\ref{sym-sp}). We have the alternative decomposition
\begin{equation}
S^2V\otimes V=S^3V\oplus\mathrm{Ker}\varphi\oplus V
\end{equation}
and therefore we can display $\Hom_{\h}(S^2V\otimes V,V)$ as
\begin{equation}\label{hom-sp}
\Hom_{\h}(S^2V\otimes V,V)=\Hom_{\h}(S^3V, V)\oplus \Hom_{\h}(\mathrm{Ker}\,\varphi,V)\oplus \Hom_{\h}(V, V)
\end{equation}
Since  $\Hom_{\h}(S^2V\otimes V,V)$ is two dimensional and $V$ is irreducible as an $\h$-module, equation (\ref{hom-sp}) shows that either $\Hom_{\h}(S^3V,V)$ is a trivial vector space or $\Hom_{\h}(S,V)$ is a one dimensional vector space spanned by $\varphi$. In case $\Hom_{\h}(S^3V,V)=0$, the triple product $xyz$ defined in (\ref{triple-sp}) restricted to $S^3V$ must be trivial. Then this product satisfies the additional identity
\begin{equation}\label{sp4}
xyz+zxy+yzx=0
\end{equation}
for all $x,y,z \in V$. Hence, using (\ref{sp1}), (\ref{sp2}) and (\ref{sp4}) we have that $(V,xyz)$ is an anti-Lie triple system with $\h$ as inner derivation algebra. Moreover, the triple system is simple and of classical type by the $\h$-irreducibility of $V$.

Otherwise, $\Hom_{\h}(S,V)=k \varphi$, and the restriction of the triple product to $S$ gives us the relationship
\begin{equation}\label{eq5.2}
xyz-zyx= \xi(2b(x,z)y+b(y,z)x-b(y,x)z)
\end{equation}
for some $\xi\in k$. Moreover $\xi$ must be nonzero: otherwise,  for all $x,y,z\in V$ we have $xyz=zyx$ and the triple products $\langle x_{\sigma}y_{-\sigma}z_{\sigma}\rangle=\sigma
x_{\sigma}y_{-\sigma}z_{\sigma}$ defined on the vector space pair $\mathcal{U}=(V^+,V^-)$ with $V^\sigma=V$ and $\sigma=\pm$, satisfy
$$
\langle x_{\sigma}y_{-\sigma}z_{\sigma}\rangle=\sigma
x_{\sigma}y_{-\sigma}z_{\sigma}=\sigma
z_{\sigma}y_{-\sigma}x_{\sigma}=\langle z_{\sigma}y_{-\sigma}x_{\sigma}\rangle
$$
and from (\ref{sp1}) and (\ref{sp2})
{\setlength\arraycolsep{2pt}
$$
\begin{array}{rcl}
\langle x_{\sigma}y_{-\sigma}\langle
u_{\sigma}v_{-\sigma}w_{\sigma}\rangle\rangle &= &\sigma^2x_{\sigma}
y_{-\sigma}(u_{\sigma}v_{-\sigma}w_{\sigma})\\
&=&\sigma^2((x_{\sigma}y_{-\sigma}u_{\sigma})v_{-\sigma}w_{\sigma})+
u_{\sigma}(y_{-\sigma}x_{\sigma}v_{-\sigma})w_{\sigma}\\
& &+ u_{\sigma}v_{-\sigma}(x_{\sigma}y_{-\sigma}w_{\sigma}))\\
&=&\langle\langle
x_{\sigma}y_{-\sigma}u_{\sigma}\rangle v_{-\sigma}w_{\sigma}\rangle-
\langle u_{\sigma}\langle
y_{-\sigma}x_{\sigma}v_{-\sigma}\rangle w_{\sigma}\rangle\\
& &+\langle u_{\sigma}v_{-\sigma}\langle
x_{\sigma}y_{-\sigma}w_{\sigma}\rangle\rangle
\end{array}
$$
}Therefore $\mathcal{U}$ is a Jordan pair for which the inner derivation operators are of the form
\begin{equation}\label{estasderbis}
\begin{array}{rl}
(D_+(x_+,y_-),D_-(y_-,x_+))&=(\langle x_+y_-\, \cdot\,\rangle, -\langle y_-x_+\, \cdot\, \rangle)\\
&=(x_+y_-\, \cdot\,,y_-x_+\, \cdot\,)
\end{array}
\end{equation}
Now from (\ref{sp1}) we have $d_{x,y}=xy\, \cdot=d_{y,x}$, thus $(D_+(x_+,y_-),D_-(y_-,x_+))=(D_{x_+,y_-},D_{x_+,y_-})$, which shows that the Lie algebra $\Inder \mathcal{U}$ is isomorphic to $\h$. Since $V$ is $\h$-irreducible, $\mathcal{U}$ is a simple Jordan pair (\cite[Proposition 1.2]{FauFe80}). But from Table \ref{simple-jordan} we deduce that the inner derivation Lie algebras of the simple Jordan pairs are not simple. Hence $\xi\neq 0$. Now, from (\ref{eq5.2}) we get $xyz-zyx=\xi(2b(x,z)y+b(y,z)x-b(y,x)z)$ for any $x,y,z$ and using (\ref{sp1}) we obtain
$$
yxz-yzx=xyz-zyx=\xi b(y,z)x-\xi b(y,x)z+2\xi b(x,z)y
$$
The previous identity, together with (\ref{sp1}), (\ref{sp2}) and (\ref{sp3}), shows that $(V,xyz)$ is a symplectic triple system with associated skew-symmetric bilinear form  $\xi b(x,y)$. Moreover, as $b(x,y)$ is nondegenerate, $V$ is a simple triple system (\cite[Proposition 2.4]{Eld06}) with $\h$ as its inner derivation algebra.
\end{proof}

Now we have all the ingredients in order to state the main result in this section:

\begin{theorem}\label{teo-sp}
Let $(\m, a\cdot b, [a,b,c])$ be an irreducible LY-algebra of generic type and standard enveloping Lie algebra of type $\frsp$. Then there is a simple symplectic triple system $(T,[...],b)$ of one of the following forms:
\begin{itemize}
\item[(i)] $\cT_k$, the symplectic triple system associated to the Jordan algebra $J=k$ with cubic form $n(\alpha)=\alpha^3$,
\item[(ii)] $\cT_{\cH_3(\cC)}$, the symplectic triple system associated to the Jordan algebra $J=\cH_3(\cC)$, where $\cC$ is either $k,k\times k, \mathrm{Mat}_2(k)$ or the algebra of octonions $\cO$.
\end{itemize}
such that, up to isomorphism, $\g(\m)=\frsp(T,b)$ and $\h=\Inder(T)$. The LY-algebra $\m$ appears as the orthogonal complement to $\h$ in $\g(\m)$ relative to the Killing form, with the binary and ternary products in \eqref{eq:binter}.
\end{theorem}

\begin{proof} Lemma \ref{estructura-sp} and Lemma \ref{casotriple-sim-aLie} show that $\m=\h^\perp$, the orthogonal complement of the inner derivation algebra $\h$ of a simple anti-Lie or symplectic triple system $(V,xyz)$ with the following extra features:
\begin{itemize}
\item[(1)] the inner derivation algebra $\Inder V$ is a simple Lie algebra,
\item[(2)]$V=V(m\lambda_i)$ as $\Inder V$-irreducible module with dominant weight $m$-times a fundamental weight $\lambda_i$ and $\dimens V \geq 4$ (actually, either $\dimens V=4$, $\h=A_1$ and $\lambda=3\lambda_1$ or $\lambda$ is a fundamental dominant weight),
\item[(3)]  $S^2V$ decomposes as a sum of two irreducible modules.
\end{itemize}
Since isomorphic irreducible anti-Lie or symplectic triple systems provide isomorphic LY-algebras, we just have to check which triple systems satisfy these extra conditions.  Following \cite{FauFe80}, the simple anti-Lie triple systems  are the odd parts of the simple Lie superalgebras and therefore the classification of such triple systems is reduced to that of the simple Lie superalgebras in \cite{Kac}. For simple symplectic triple systems we will follow the classification and comments given in \cite[Section 2]{Eld06}. Both classifications are outlined in the Appendix of the paper. Table \ref{simple-superalgebras} shows that there are no simple  anti-Lie triple systems satisfying (1)--(3) simultaneously. For simple symplectic triple systems, following Table \ref{simple-simplecticos} and applying  the restrictions (1)-(2)-(3), we get that the only possibilities are given by the simple symplectic triple systems $\cT_k$, associated to the one dimensional Jordan algebra $k$ in item (i), and the simple symplectic triple systems $\cT_\cJ$ associated to a simple Jordan algebra $\cJ=\cH_3(\cC)$ of degree $3$ with $\cC=k, k\times k$, the algebra of quaternions or the algebra of octonions, which proves the theorem.
\end{proof}


\section{Exceptional case}

In this section we deal with the irreducible LY-algebras of Generic Type and exceptional standard enveloping Lie algebra. These systems appear in reductive decompositions $\g=\h\oplus \m$  for  which (a), (b) and (c) in (\ref{condicionescasonosimpleI}) hold,  $\g$ being a simple Lie algebra of type $G_2$, $F_4$, $E_6$, $E_7$, or $E_8$. The classification in this case is given in the next result

\begin{theorem}
 Let $(\m, x\cdot y,[x,y,z])$ be an irreducible LY-algebra of generic type and  exceptional standard enveloping algebra. Then one of the following holds:
\begin{itemize}
\item[(i)] $\m$ is the Lie triple system associated to one of the symmetric pairs $(F_4, B_4)$, $(E_{6}, F_4)$, $(E_6, C_4)$, $(E_7, A_7)$ or $(E_8,D_8)$.
\item[(ii)] $\m=\h^\perp$ is the orthogonal complement with respect to the Killing form in $\g$ associated to one of the reductive pairs $(\g,\h)=(G_2,A_1)$, $(E_6,G_2)$ or $(E_7,A_2)$. Moreover, for the previous pairs, as a module for $\h$, $\m$ is isomorphic to $V(10\lambda_1)$, $V(\lambda_1\oplus \lambda_2)$ and $V(4\lambda_1\oplus 4\lambda_2)$ respectively.
\end{itemize}
\end{theorem}

\begin{proof}

The systems appear in reductive decompositions $\g=\h\oplus \m$  satisfying  (\ref{condicionescasonosimpleI}). Lemma 4.2 in \cite{BEM09} shows that for each such reductive decomposition there exists an analogous decomposition $\tilde\g=\tilde\h\oplus\tilde\m$ over the complex numbers. In particular, the highest weight of $\m$ as a module for $\h$ coincides with the highest weight of $\tilde\m$ as a module for $\tilde\h$. Then, because of \cite[Lemma 4.3]{BEM09}, either $\m$ is a Lie triple system and we obtain (i), or $\h$ is a simple $S$-subalgebra of $\g$ and the different possibilities for the pair $(\g,\h)$ can be read from \cite[Theorem 14.1]{Dyn}, where a complete list of the complex maximal and simple $S$-subalgebras of the exceptional Lie algebras is given: $(\g,\h)=(G_2,A_1)$, $(F_4,A_1)$, $(E_6,A_1)$, $(E_7,A_1)$, $(E_8,A_1)$, $(E_6,G_2)$, $(E_6,C_4)$, $(E_6,F_4)$, or $(E_7,A_2)$.

The cases $(E_6,C_4)$ and $(E_6,F_4)$ correspond to symmetric pairs already considered in (i). In case $\h=A_1$, the irreducibility restriction on $\m$ forces $\m=V(n\lambda_1)$ for some $n\geq 1$. Now, given a Cartan subalgebra, spanned by an element $h$, of $A_1$ we can pick a Cartan subalgebra $H$ of $\g$ with $h\in H$. Then, $[H,h]=0$, so $H \subseteq C_{\g}(h)=k\cdot h\oplus V(n\lambda_1)_0$, where $V(n\lambda_1)_0$ is the $0$-weight subspace of $V(n\lambda_1)$. Since $\dimens V(n\lambda_1)_0$ is $0$ or $1$, depending on the parity of $n$, we get that $\g$ must be a Lie algebra of rank $2$. Hence, $(G_2,A_1)$ is the unique possibility that works. A dimension count shows that $n=10$ in this case.

A dimension count for the other cases $(E_6,G_2)$ and $(E_7,A_2)$ completes the proof.
\end{proof}

\begin{remark}
The reductive pair $(G_2,A_1)$ can be constructed by using transvections. A construction of $G_2$ from $A_1=\spl_2(k)$ and an eleven-dimensional module is given in Dixmier \cite{Dixmier} (see also \cite{BDE}).
The symmetric pairs $(F_4,B_4)$ and $(E_6,F_4)$ are strongly related to the Albert algebra (the exceptional simple Jordan algebra, see \cite{JacJA}). Nice constructions of the pairs $(E_8,D_8)$, $(E_7,A_7)$ and $(E_6,C_4)$ can be read from constructions in \cite{Adams}.
\end{remark}

\clearpage

\section{Appendix}


\subsection{Lie and anti-Lie triple systems}\label{Lie} Following \cite[Theorem 1.1]{Lis} Lie triple systems are nothing else but skew-symmetric elements relative to involutive automorphisms in Lie algebras. That is, these systems can be viewed as the odd part of $\Z_2$-graded Lie algebras. Anti-Lie triple systems appear in the same vein by using the odd part of Lie superalgebras. Since Lie triple systems are LY-algebras with trivial binary product, because of Definition \ref{df:LY} of this paper and \cite[Section 5]{FauFe80}, it is possible to introduce both triple systems in an axiomatically unified way by using a vector space $V$ endowed with a triple product $xyz$ satisfying the identities
\begin{equation}\label{lie-antilie}
\begin{array}{l}
\quad xyz=\epsilon\, yxz\\
\quad xyz+yzx+zxy=0\\
\quad  xy(uvw)= (xyu)vw+u(xyv)w+uv(xyw)\\
\end{array}
\end{equation}
where $\epsilon=-1$  for Lie triple systems and $\epsilon=1$ for anti-Lie triple systems.

Given a Lie or anti-Lie triple system $(V,xyz)$, the standard enveloping construction $\g(V)=D(V,V)\oplus V$ in (\ref{eq:gm}), where $D(V,V)=\Inder V=\eespan \langle xy\, \cdot: x,y \in V\rangle$ is the inner derivation Lie algebra of the corresponding triple system, provides either a  $\Z_2$-graded Lie algebra or a Lie superalgebra according to $V$ being a Lie or and anti-Lie triple system. Moreover, it is not difficult to prove that the Lie algebra (respectively superalgebra) $\g(V)$ is graded simple (respectively simple) if and only if the Lie (respectively anti-Lie) triple system $V$ is simple.

Over algebraically closed fields of characteristic zero, simple Lie triple systems were classified in \cite{Lis} through involutive automorphisms. Table I in \cite{Fau80} presents an alternative classification of these systems by means of (reductive) symmetric pairs $(\g(V), \Inder V)$ obtained from affine Dynkin diagrams (see \cite[Chapter 4]{Kac90}) that encode the Cartan type of the standard enveloping Lie algebra $\g(V)$ and the inner derivation Lie algebra $\Inder V$ of the Lie triple system $V$. These diagrams are equipped with some numerical labels which describe the lowest  weight of $V$ as $\Inder V$-module (the highest weight also is easily checked since simple Lie triple systems are self dual modules). Using the latter classification, in
\cite[Table III]{Fau85} all simple and $\Inder V$-irreducible Lie triple systems are listed. Combining the results in \cite{Fau80} and \cite{Fau85} we arrive at Table \ref{ltsadjuntos-so-sl-excep}, that displays the irreducible Lie triple systems for which its inner derivation algebra is simple.

On the other hand, as simple anti-Lie triple systems are the odd part of simple Lie superalgebras, the classification of these systems can be obtained from that of the simple Lie superalgebras in \cite{Kac}. Of special interest for our purposes are the simple anti-Lie triple systems $V$ which are completely reducible as modules for $\Inder V$. This class of systems appears from the odd parts of the simple classical Lie superalgebras listed in Table \ref{simple-superalgebras}. We shall refer to them as anti-Lie triple systems of {\em classical type}. Structural module information on Table \ref{simple-superalgebras} follows from \cite[Propositions 2.1.2]{Kac}.

\subsection{Symplectic and orthogonal triple systems}\label{simplecticos}

Symplectic triple systems were introduced in \cite{YamAs} and orthogonal triple systems were defined in \cite[Section V]{Oku93}. They are basic ingredients in the construction of some $5$-graded Lie algebras and Lie superalgebras respectively (see \cite{Eld06}), and hence they are strongly related to $\Z_2$-graded Lie algebras and to a specific class of Lie superalgebras. These triple systems consist of a vector space $V$ endowed with a trilinear product $xyz$ and a $\epsilon$-bilinear form $b$, with $\epsilon=-1$ (skewsymmetric)  for symplectic triple systems and $\epsilon=1$ (symmetric) for orthogonal ones, satisfying the relations
\begin{equation}\label{sym-ort}
\begin{array}{l}
\quad xyz=-\epsilon\, yxz\\
\quad xyz+\epsilon xzy=\epsilon\,b(x,y)z+b(x,z)y-\epsilon\,2b(y,z)x\\
\quad  xy(uvw)= (xyu)vw+u(xyv)w+uv(xyw)\\
\end{array}
\end{equation}
We note that the second relation in the orthogonal case ($\epsilon=1$) is just the linearization of the identity
\begin{equation}
\quad xyy=b(x,y)y- (y,y)x
\end{equation}
and, from the third relation, we can introduce for these systems in the usual way the inner derivation Lie algebra $\Inder V=\eespan \langle xy\, \cdot: x,y \in V\rangle$.

Symplectic and orthogonal triple systems are related to the so called $(-\epsilon,-\epsilon)$ balanced Freudenthal-Kantor triple systems introduced in \cite{YamOn}. In \cite[Theorems 2.16 and 2.18]{Eld06} it is also shown that symplectic triple systems are closely related to Freudenthal triple systems and  a class of ternary algebras defined in \cite{FauFe}: the balanced symplectic Lie algebras. Moreover, following \cite[Theorems 2.4 and 4.4]{Eld06}, the simplicity of both types of triple systems (fields of characteristic different from $2$ and $3$) is equivalent to the non-degeneracy of the associated bilinear form $b$.

The relationship between symplectic triple systems, Freudenthal triple systems and ternary algebras leads to the classification of simple symplectic triple systems over algebraically closed fields of characteristic different from $2$ and $3$ given in \cite[Theorem 2.21]{Eld06} (the classification is based on previous classifications of Freudenthal triple systems in \cite{Mey68} and simple ternary algebras from \cite[Therorem 4.1]{FauFe}). For simple orthogonal triple systems, Theorem 4.7 in \cite{Eld06} displays the classification over algebraically closed fields of characteristic zero, by means of a previous classification of the simple $(-1,-1)$ balanced Freudenthal-Kantor triple systems in \cite[Theorem 4.3]{EldKaO}. The classifications and comments therein \cite{Eld06} provide Tables \ref{simple-simplecticos} and \ref{simple-ortogonales}, where the Cartan type of the Lie algebra $\Inder V$ and the $\Inder V$-module structure of $V$ is given for the different types of simple symplectic and orthogonal triple systems.

Among the simple symplectic triple systems a special use
of the following ones will be made:
\begin{equation}\label{eq:TJsymplectic}
\cT_{\J}=\Bigl\{\begin{pmatrix} \alpha & a \\ b &
 \beta \end{pmatrix}: \alpha, \beta \in k, a, b \in \J\Bigr\}
\end{equation}
where $\J=\J ordan(n,c)$ is the Jordan algebra of a nondegenerate
cubic form $n$ with basepoint (see \cite[II.4.3]{McC04} for a
definition) of one of the following types: $\J=k,
n(\alpha)=\alpha^3$ and $t(\alpha, \beta)=3\alpha \beta$ or
$\J=H_3(\cC)$ for a unital composition algebra $\cC$. Theorem 2.21
in \cite{Eld06} displays the product and bilinear form for
the triple systems $\cT_\J$ by using the trace form $t(a,b)$ and the cross product $a\times b$ attached to the Jordan algebra $\cJ$.

\subsection{Jordan and anti-Jordan pairs}\label{pares} Jordan pairs, axiomatically introduced (for arbitrary fields and dimension) in \cite{Loos} are basic ingredients in the construction of Lie algebras with short $3$-gradings. In the context of Lie superalgebras endowed with a consistent short $3$-grading, the anti-Jordan pairs introduced in \cite{FauFe80} constitute the corresponding concept. Following \cite{FauFe80} (see also \cite[Chapter XI]{Mey72}), over fields of characteristic different from $2$ and $3$ both types of pairs can be defined by means of a pair of vector spaces $\mathcal{U}=(\mathcal{U}^+,\mathcal{U}^-)$ with trilinear products $\{x_{\sigma}y_{-\sigma}z_{\sigma}\}$ for $\sigma=\pm$ satisfying the identities:

\begin{equation}\label{jp-ajp}
\begin{array}{l}
\quad \{x_{\sigma}y_{-\sigma}z_{\sigma}\}=\epsilon \{z_{\sigma}y_{-\sigma}x_{\sigma}\}\\
\quad  \{x_{\sigma}y_{-\sigma}\{u_{\sigma}v_{-\sigma}w_{\sigma}\}\}= \{\{x_{\sigma}y_{-\sigma}u_{\sigma}\}v_{-\sigma}w_{\sigma}\}-\\
\qquad \qquad \qquad \qquad \qquad  \epsilon\{u_{\sigma}\{y_{-\sigma}x_{\sigma}v_{-\sigma}\}w_{\sigma}\}+ \{u_{\sigma}v_{-\sigma}\{x_{\sigma}y_{-\sigma}w_{\sigma}\}\}
\end{array}
\end{equation}
with $\epsilon=1$ for Jordan pairs and $\epsilon=-1$ for anti-Jordan pairs.

The simplicity of these systems can be characterized through its inner derivation algebra:
\begin{equation}\label{inder-jp-ajp}
\Inder\, \mathcal{U}=\eespan \langle(\{x_+y_-\cdot\},-\epsilon\{y_-x_+\cdot\}): x_+\in \mathcal{U}^+, y_-\in \mathcal{U}^-\rangle
\end{equation}
which is a Lie subalgebra of $\gl(\mathcal{U}^+)\times \gl(\mathcal{U}^-)$. According to \cite[Proposition 1.2]{FauFe80}, $\mathcal{U}$ is a simple pair if and only if $\{\mathcal{U}^\sigma \mathcal{U}^{-\sigma}\mathcal{U}^\sigma\}\neq 0$ and $\mathcal{U}^\sigma$ is an irreducible $\Inder\, \mathcal{U}$-module (via the action of the $\sigma$-component). Over algebraically closed fields of characteristic zero, the simple finite-dimensional Jordan and anti-Jordan pairs were classified in \cite[Theorem 17.12]{Loos} and
\cite[Sections 3 and 4]{FauFe80}. In Tables \ref{simple-jordan} and \ref{simple-antijordan} below a complete description of both classifications is given. The tables include the inner derivation algebra $\Inder\, \mathcal{U}$, the Cartan type of its derived subalgebra $\Inder_{0}\, \mathcal{U}=[\Inder\, \mathcal{U},\Inder\, \mathcal{U}]$ and the highest weight of $\mathcal{U}^\sigma$ as a module for $(\Inder\, \mathcal{U})'$ for the different simple Jordan and anti-Jordan pairs $\mathcal{U}$. We follow the matricial description of the original classifications, although alternative descriptions could be displayed. In this way, rectangular $p\times q$ matrices, $n\times n$ symmetric or alternating matrices are represented in the tables by $\cM_{p,q}(k)$, $\cH_n(k)$ and $\cA_n(k)$, while $\cM_{1,2}(\cO)$ represents the space of $1\times 2$ matrices over the octonions $\cO$. We use the standard notation $\cH_3(\cO)$ for the 27-dimensional exceptional Jordan algebra (or Albert algebra, see  \cite{JacJA} or \cite{Sch} for a complete description). For a given matrix $y$, its transpose is denoted by $y^t$ and in case $y\in \cM_{1,2}(\cO)$, $\bar{y}$ represents the standard involution induced in $\cM_{1,2}(\cO)$ by the involution of $\cO$. Triple products for Jordan pairs and anti-Jordan pairs of the form $\mathcal{U}=(k^n,k^n)$ are defined by means of the operators $b_{x,y}=b(y,\cdot)x-\epsilon b(x,\cdot)y$ for a nondegenerate $\epsilon$-symmetric form $b$ ($\epsilon =1$  for $b$ symmetric and $\epsilon=-1$ in case $b$ is skewsymmetric).

The structural module information given on these tables can be obtained from a direct computation for the different Jordan and anti-Jordan pairs. Alternatively, the relationship among $\Z_2$-graded simple Lie algebras (respectively superalgebras) having a consistent short $3$-grading and Jordan pairs (resp. anti-Jordan pairs) allows us to obtain the complete information from the classification  of simple and non irreducible Lie triple systems (using the corresponding affine Dynkin diagrams in \cite[Table I]{Fau80}) and simple superalgebras of type $A(m,n), C(n)$ and $P(n)$ (from \cite[Proposition 2.1.2]{Kac})

\begin{table}[htb]
\caption{Irreducible L.t.s. with simple inner derivation Lie algebra}\label{ltsadjuntos-so-sl-excep}
\medskip
\hskip -0.5cm
\begin{minipage}{\linewidth}
\begin{tabular}{lllll}
\hline

& & & &\\[-2ex]
{\scriptsize $ (\g(V), \Inder\, V)$} & {\scriptsize $\Inder\, V$} & {\scriptsize  $V=V(k\lambda_i)$} &
{\scriptsize  $V(2k\lambda_i-\alpha_i)$} & {\scriptsize  $\so(V) \cong \wedge^2 V$}\\
& & & &\\[-2ex]
\hline
&&&&\\[-2ex]

 {\scriptsize  $(A_1\times A_1,A_1)$} & {\scriptsize  $V(2\lambda_1)$} & {\scriptsize  $V(2\lambda_1)$} &
{\scriptsize $V(2\lambda_1)$} & {\scriptsize $V(2\lambda_1)$}\\

{\scriptsize  $(A_n\times A_n,A_n)_{n\geq 2}$}
& {\scriptsize  $V(\lambda_1\tiny{+}\lambda_n)$} & {\scriptsize  $V(\lambda_1\tiny{+}\lambda_n)$} &
  & \\

 {\scriptsize  $(B_3\times B_3,B_3)$} & {\scriptsize  $V(\lambda_2)$} & {\scriptsize  $V(\lambda_2)$} &
{\scriptsize $V(\lambda_1+2\lambda_3)$} & {\scriptsize $V(\lambda_2)\oplus V(\lambda_1+2\lambda_3)$}\\

{\scriptsize  $(B_n\times B_n,B_n)_{n\geq 4}$} & {\scriptsize  $V(\lambda_2)$} & {\scriptsize  $V(\lambda_2)$} &{\scriptsize $V(\lambda_1+\lambda_3)$}& {\scriptsize $V(\lambda_2)\oplus V(\lambda_1+2\lambda_3)$}\\

{\scriptsize  $(C_n\times C_n,C_n)_{n\ge 2}$} & {\scriptsize  $V(2\lambda_1)$} & {\scriptsize  $V(2\lambda_1)$} &
{\scriptsize $V(2\lambda_1+\lambda_2)$} & {\scriptsize $V(2\lambda_1)\oplus V(2\lambda_1+\lambda_2)$}\\

 {\scriptsize  $(D_4\times D_4,D_4)$} & {\scriptsize  $V(\lambda_2)$} & {\scriptsize  $V(\lambda_2)$} &
{\scriptsize $V(\lambda_1\tiny{+}\lambda_3\tiny{+}\lambda_4)$} & {\scriptsize $V(\lambda_2)\oplus V(\lambda_1\tiny{+}\lambda_3\tiny{+}\lambda_4)$}\\

{\scriptsize  $(D_n\times D_n,D_n)_{n\geq 5}$} & {\scriptsize  $V(\lambda_2)$} & {\scriptsize  $V(\lambda_2)$} &
{\scriptsize $V(\lambda_1+\lambda_3)$} & {\scriptsize $V(\lambda_2)\oplus V(\lambda_1+2\lambda_3)$}\\

{\scriptsize  $(G_2\times G_2,G_2)$} & {\scriptsize  $V(\lambda_2)$} & {\scriptsize  $V(\lambda_2) $} &
{\scriptsize  $V(3\lambda_1)$} &
{\scriptsize  $V(\lambda_2)\oplus V(3\lambda_1)$}\\

{\scriptsize  $(F_4\times F_4,F_4)$} & {\scriptsize  $V(\lambda_1)$} & {\scriptsize  $V(\lambda_1) $} &
{\scriptsize  $V(\lambda_2)$} &
{\scriptsize  $V(\lambda_1)\oplus V(\lambda_2)$}\\

{\scriptsize  $(E_6\times E_6,E_6)$} & {\scriptsize  $V(\lambda_2)$} & {\scriptsize  $V(\lambda_2)$} &
{\scriptsize  $V(\lambda_4)$} &
{\scriptsize  $V(\lambda_2)\oplus V(\lambda_4)$} \\

{\scriptsize  $(E_7\times E_7,E_7)$} & {\scriptsize  $V(\lambda_1)$} & {\scriptsize  $V(\lambda_1)$} &
{\scriptsize  $V(\lambda_3)$} &
{\scriptsize  $V(\lambda_1)\oplus V(\lambda_3)$} \\

{\scriptsize  $(E_8\times E_8,E_8)$} & {\scriptsize  $ V(\lambda_8)$} & {\scriptsize  $V(\lambda_8)$} &
{\scriptsize  $V(\lambda_7) $} &
{\scriptsize  $V(\lambda_8)\oplus V(\lambda_7)$}\\
& & & &\\[-2ex]
& & & &\\[-2ex]

{\scriptsize  $(D_3,B_2)$} & {\scriptsize  $V(2\lambda_2)$} &
{\scriptsize  $V(\lambda_1)$} &
{\scriptsize  $V(2\lambda_2)$} &
{\scriptsize  $V(2\lambda_2)$}\\

 {\scriptsize  $(B_3,D_3)$} & {\scriptsize  $V(\lambda_2\tiny{+}\lambda_3)$} & {\scriptsize  $V(\lambda_1)$} &
{\scriptsize $V(\lambda_2+\lambda_3)$} & {\scriptsize $V(\lambda_2+\lambda_3)$}\\

 {\scriptsize  $(B_n,D_n)_{n\geq 4}$} & {\scriptsize  $V(\lambda_2)$} & {\scriptsize  $V(\lambda_1)$} &
{\scriptsize $V(\lambda_2)$} & {\scriptsize $V(\lambda_2)$}\\

 {\scriptsize  $(D_{n+1},B_n)_{n\geq 3}$} & {\scriptsize  $V(\lambda_2)$} & {\scriptsize  $V(\lambda_1)$} &
{\scriptsize $V(\lambda_2)$} & {\scriptsize $V(\lambda_2)$}\\
& & & &\\[-2ex]
& & & &\\[-2ex]

{\scriptsize  $(A_2,A_1)$} & {\scriptsize  $V(2\lambda_1)$} & {\scriptsize  $V(4\lambda_1)$} &
{\scriptsize $V(6\lambda_1)$} & {\scriptsize $V(2\lambda_1)\oplus V(6\lambda_1)$}\\

{\scriptsize  $(A_4,B_2)$} & {\scriptsize  $V(2\lambda_2)$} & {\scriptsize  $V(2\lambda_1)$} &
{\scriptsize $V(2\lambda_1+2\lambda_2)$} & {\scriptsize $V(2\lambda_2)\oplus V(2\lambda_1+2\lambda_2)$}\\

{\scriptsize  $(A_{2n},B_n)_{n\geq 3}$} & {\scriptsize  $V(\lambda_2)$} & {\scriptsize  $V(2\lambda_1)$} &
{\scriptsize $V(2\lambda_1+\lambda_2)$} & {\scriptsize $V(\lambda_2)\oplus V(2\lambda_1+\lambda_2)$}\\

{\scriptsize  $(A_{2n-1},C_n)_{n\geq 3}$} & {\scriptsize  $V(2\lambda_1)$} & {\scriptsize  $V(\lambda_2)$} &
{\scriptsize $V(\lambda_1+\lambda_3)$} & {\scriptsize $V(2\lambda_1)\oplus V(\lambda_1+\lambda_3)$}\\

{\scriptsize  $(A_5,D_3)_{n\geq 4}$} & {\scriptsize  $V(\lambda_2\tiny{+}\lambda_3)$} & {\scriptsize  $V(2\lambda_1)$} &
{\scriptsize $V(2\lambda_1\tiny{+}\lambda_2\tiny{+}\lambda_3)$} & {\scriptsize $V(\lambda_2\tiny{+}\lambda_3)\oplus V(2\lambda_1\tiny{+}\lambda_2\tiny{+}\lambda_3)$}\\

{\scriptsize  $(A_{2n-1},D_n)_{n\geq 4}$} & {\scriptsize  $V(\lambda_2)$} & {\scriptsize  $V(2\lambda_1)$} &
{\scriptsize $V(2\lambda_1+\lambda_2)$} & {\scriptsize $V(\lambda_2)\oplus V(2\lambda_1+\lambda_2)$}\\

& & & &\\[-2ex]
& & & &\\[-2ex]

{\scriptsize  $(E_6,F_4)$} & {\scriptsize  $V(\lambda_1)$} & {\scriptsize  $V(\lambda_4)$} &
{\scriptsize $V(\lambda_3)$} & {\scriptsize $V(\lambda_1)\oplus V(\lambda_3)$}\\

{\scriptsize  $(F_4,B_4)$} & {\scriptsize  $V(\lambda_2)$} & {\scriptsize  $V(\lambda_4)$} &
{\scriptsize $V(\lambda_3)$} & {\scriptsize $V(\lambda_2)\oplus V(\lambda_3)$}\\

{\scriptsize  $(E_6,C_4)$} & {\scriptsize  $V(2\lambda_1)$} & {\scriptsize  $V(\lambda_4)$} &
{\scriptsize $V(2\lambda_3)$} & {\scriptsize $V(2\lambda_1)\oplus V(2\lambda_3)$}\\

{\scriptsize  $(E_7,A_7)$} & {\scriptsize  $V(\lambda_1\tiny{+}\lambda_7)$} & {\scriptsize  $V(\lambda_4)$} &
{\scriptsize $V(\lambda_3+\lambda_5)$} & {\scriptsize $V(\lambda_1+\lambda_7)\oplus V(\lambda_3+\lambda_5)$}\\

{\scriptsize  $(E_8,D_8)$} & {\scriptsize  $V(\lambda_2)$} & {\scriptsize  $V(\lambda_8)$} &
{\scriptsize $V(\lambda_6)$} & {\scriptsize $V(\lambda_2)\oplus V(\lambda_6)$}\\
& & & &\\[-2ex]
\hline
\end{tabular}
\end{minipage}
\end{table}

\renewcommand{\thempfootnote}{\fnsymbol{mpfootnote}}

\begin{table}[htb]
\vskip 1cm
\caption{Simple classical Lie superalgebras}\label{simple-superalgebras}
\medskip
\begin{minipage}{\linewidth}
\begin{center}
\begin{tabular}{lll}
\hline
& &  \\[-2ex]
{\scriptsize $\cL$-Type} &{\scriptsize  $\cL_{\bar{0}}$} &
{\scriptsize  $\cL_{\bar{1}}$ as $\cL_{\bar{0}}$-module}\\

& &  \\[-2ex]
\hline
&&\\[-2ex]

 {\scriptsize  $A(m,0)_{m\ge 1}$} & {\scriptsize  $A_m\times Z$}\footnote{\ $Z$ stands for a one-dimensional center of $\cL_{\bar 0}$. In case $\cL_{\bar 0}=[\cL_{\bar 0},\cL_{\bar 0}]\times Z$, the highest weight of $\cL_{\bar 1}$ as a module for $[\cL_{\bar 0},\cL_{\bar 0}]$ is considered. For these cases, $\cL_{\bar 1}$ decomposes as sum of two irreducible modules. The elements of the center act as $\alpha \cdot \mathrm{Id}$ in one of the two summands, $\alpha$ being a nonzero scalar, and  as $-\alpha \cdot \mathrm{Id}$ on the other summand. The same remark works for the remaining tables.} &
{\scriptsize $V(\lambda_1)\oplus V(\lambda_{m})$}\\
 {\scriptsize  $A(m,n)_{m>n\ge 1}$} & {\scriptsize  $A_m\times A_n\times Z$} &
{\scriptsize $V(\lambda_1)\otimes V(\lambda_{1}')\oplus V(\lambda_{m})\otimes V(\lambda_{n}')$}\\
& & \\[-2ex]
 {\scriptsize  $A(n,n)_{n\ge 1}$}& {\scriptsize  $A_{n}\times A_n$} &
{\scriptsize $V(\lambda_1)\otimes V(\lambda_{1}')\oplus V(\lambda_{n})\otimes V(\lambda_{n}')$}\\
& & \\[-2ex]

 {\scriptsize  $B(0,n)_{n\ge 1}$}& {\scriptsize $C_n$}& {\scriptsize $ V(\lambda_1)$}\\


{\scriptsize  $B(m,n)_{m, n\ge 1}$}& {\scriptsize $B_m\times C_n$}& {\scriptsize $ V(\lambda_1)\otimes V(\lambda_1')$}\\

& &\\[-2ex]

 {\scriptsize  $D(2,n)_{n\ge 2}$}& {\scriptsize $A_1\times A_1\times C_n$}&
{\scriptsize $ V(\lambda_1)\otimes V(\lambda_1')\otimes V(\lambda_1'')$}\\

 {\scriptsize  $D(m,n)_{m\ge 3, n\ge 1}$}& {\scriptsize $D_m\times C_n$}&
{\scriptsize $ V(\lambda_1)\otimes V(\lambda_1')$}\\

& &\\[-2ex]

{\scriptsize  $C(n)_{n\ge 2}$} &
{\scriptsize $C_{n-1}\times Z$}&{\scriptsize $V(\lambda_1)\oplus V(\lambda_1)$}\\

& &\\[-2ex]

{\scriptsize  $Q(n)_{n\ge 2}$} &
{\scriptsize $A_{n}$}&{\scriptsize $V(\lambda_1+\lambda_n)$}\\

& &\\[-2ex]

{\scriptsize  $P(n)_{n\ge 2}$}&{\scriptsize  $A_n$} &
{\scriptsize $V(\lambda_{n-1})\oplus V(2\lambda_1)$}\\
& &\\[-2ex]

 {\scriptsize  $D(2,1;\alpha)_{\alpha\neq 0,-1}$} &
{\scriptsize $A_1\times A_1\times A_1$}&{\scriptsize $V(\lambda_1)\otimes V(\lambda_1')\otimes V(\lambda_1'')$}\\

& &\\[-2ex]

{\scriptsize  $F(4)$} & {\scriptsize  $B_3\times A_1$} &
{\scriptsize $V(\lambda_3)\otimes V(\lambda_1')$}\\

 & &\\[-2ex]
{\scriptsize  $G(3)$}& {\scriptsize  $G_2\times A_1$} &
{\scriptsize $V(\lambda_1)\otimes V(\lambda_1')$}\\

& & \\[-2ex]

\hline

\end{tabular}

\end{center}
\end{minipage}

\end{table}

\begin{table}[htb]
\vskip 1cm
\caption{Simple symplectic triple systems}\label{simple-simplecticos}
\medskip
\begin{minipage}{\linewidth}
\begin{center}
\begin{tabular}{lllll}
\hline
& & & \\[-2ex]
{\scriptsize $V$-Type} &{\scriptsize  $\dimens V$} & {\scriptsize  $\Inder V$} &
{\scriptsize  $V$ as $\Inder V$-module}& {\scriptsize  $\frsp(V)\cong S^2V$}\\

& & & &\\[-2ex]
\hline
&&&&\\[-2ex]

 {\scriptsize  Orthogonal-type} & {\scriptsize  $8$} & {\scriptsize  $A_{1}\tiny{\times}A_1\tiny{\times} A_1$} &
{\scriptsize $V(\lambda_1)\tiny{\otimes} V(\lambda_{1}')\tiny{\otimes} V(\lambda_{1}'')$}&\\[-0,5ex]
 & {\scriptsize  $4n, n \geq 3$} & {\scriptsize  $A_{1}\tiny{\times} D_n$} &
{\scriptsize $V(\lambda_1)\tiny{\otimes} V(\lambda_{1}')$}&\\[-0,5ex]

&{\scriptsize $6$} & {\scriptsize $A_1\tiny{\times} A_1$}& {\scriptsize $ V(\lambda_1)\tiny{\otimes} V(2\lambda_1')$}& \\[-0,5ex]
 &{\scriptsize $4n+2, n \geq 2$} &{\scriptsize  $A_1\tiny{\times} B_n$} &
{\scriptsize $V(\lambda_1)\tiny{\otimes}V(\lambda_{1}')$}&{}\\[-0,5ex]

& & & &\\[-1ex]
{\scriptsize  Special-type} &{\scriptsize $2$} &{\scriptsize  $ Z$} &
{\scriptsize $k\times k$}&\\[-0,5ex]
&{\scriptsize $2n, n\geq 2$} &{\scriptsize  $A_{n-1}\times Z$} &
{\scriptsize $V(\lambda_1)\otimes V(\lambda_{n-1})$}&\\

& & & &\\[-1ex]

 {\scriptsize   Symplectic-type} & {\scriptsize  $2$} & {\scriptsize  $A_{1}$} &
{\scriptsize $V(\lambda_1)$}&{\scriptsize $V(2\lambda_1)$}\\ [-0,5ex]
& {\scriptsize  $2n, n \geq 2$} & {\scriptsize  $C_n$} &
{\scriptsize $V(\lambda_1)$}&{\scriptsize $V(2\lambda_1)$}\\

& & & &\\[-1ex]

 {\scriptsize   $\cT_k$} & {\scriptsize  $4$} & {\scriptsize  $A_{1}$} &
{\scriptsize $V(3\lambda_1)$}&{\scriptsize $V(6\lambda_1)\oplus V(2\lambda_1)$}\\


{\scriptsize   $\cT_{\cH_3(k)}$} & {\scriptsize  $14$} & {\scriptsize  $C_3$} &
{\scriptsize $V(\lambda_3)$}&{\scriptsize $V(2\lambda_1)\oplus V(2\lambda_3)$}\\

{\scriptsize   $\cT_{\cH_3(k\times k)}$} & {\scriptsize  $20$} & {\scriptsize  $A_5$} &
{\scriptsize $V(\lambda_3)$}&{\scriptsize $V(\lambda_1+\lambda_5)\oplus V(2\lambda_3)$}\\

{\scriptsize   $\cT_{\cH_3(\cQ)}$} & {\scriptsize  $32$} & {\scriptsize  $D_6$} &
{\scriptsize $V(\lambda_6)$}&{\scriptsize $V(\lambda_2)\oplus V(2\lambda_6)$}\\


{\scriptsize   $\cT_{\cH_3(\cO)}$} & {\scriptsize  $56$} & {\scriptsize  $E_7$} &
{\scriptsize $V(\lambda_7)$}&{\scriptsize $V(\lambda_1)\oplus V(2\lambda_7)$}\\

& & & &\\[-1ex]

\hline
\end{tabular}
\end{center}
\end{minipage}
\end{table}

\begin{table}[htb]
\vskip 1cm
\caption{Simple orthogonal triple systems}\label{simple-ortogonales}
\medskip
\begin{minipage}{\linewidth}
\begin{center}
\begin{tabular}{lllll}
\hline
& & & \\[-2ex]
{\scriptsize $V$-Type} &{\scriptsize  $\dimens V$} & {\scriptsize  $\Inder V$} &
{\scriptsize  $V$ as $\Inder V$-module}& {\scriptsize  $\so(V)\cong \wedge^2V$}\\

& & & &\\[-2ex]
\hline
&&&&\\[-2ex]

 {\scriptsize  Orthogonal-type} &{\scriptsize $3,5$} &{\scriptsize  $A_1,B_2$} &
{\scriptsize $V(2\lambda_1),V(\lambda_1)$}&{\scriptsize   $V(2\lambda_1), V(2\lambda_2)$}\\[-0,5ex]
 &{\scriptsize $2n+1, n \geq 3$} &{\scriptsize  $B_n$} &
{\scriptsize $V(\lambda_1)$}&{\scriptsize   $V(\lambda_2)$}\\[-0,5ex]
&{\scriptsize $2,4$} & {\scriptsize $Z, A_1\times A_1$}& {\scriptsize $k\times k, V(\lambda_1)\otimes V(\lambda_1')$}& \\[-0,5ex]
&{\scriptsize $6$} & {\scriptsize $D_3$}&{\scriptsize $V(\lambda_1)$}& {\scriptsize  $V(\lambda_2+\lambda_3)$}\\[-0,5ex]
&{\scriptsize $2n, n \geq 4$} & {\scriptsize $D_n$}&{\scriptsize $V(\lambda_1)$}& {\scriptsize  $V(\lambda_2)$}\\
& & & &\\[-1ex]

 {\scriptsize  Unitarian-type} & {\scriptsize  $2n, n \geq 3$} & {\scriptsize  $A_{n-1}\times Z$} &
{\scriptsize $V(\lambda_1)\oplus V(\lambda_{n-1})$}&\\
& & & &\\[-1ex]

 {\scriptsize   Symplectic-type} & {\scriptsize  $4n, n \geq 2$} & {\scriptsize  $A_{1}\times C_n$} &
{\scriptsize $V(\lambda_1)\otimes V(\lambda_1')$}&\\

& & & &\\[-1ex]

 {\scriptsize   $D_\mu$-type} & {\scriptsize  $4$} & {\scriptsize  $A_{1}\times A_1$} &
{\scriptsize $V(\lambda_1)\oplus V(\lambda_1')$}&\\

& & & &\\[-1ex]

{\scriptsize   $G$-type} & {\scriptsize  $7$} & {\scriptsize  $G_2$} &
{\scriptsize $V(\lambda_1)$}&{\scriptsize $V(\lambda_1)\oplus V(\lambda_2)$}\\

& & & &\\[-1ex]

{\scriptsize   $F$-type} & {\scriptsize  $8$} & {\scriptsize  $B_3$} &
{\scriptsize $V(\lambda_3)$}&{\scriptsize $V(\lambda_1)\oplus V(\lambda_2)$}\\

& & & &\\[-1ex]

\hline
\end{tabular}
\end{center}
\end{minipage}
\end{table}

\begin{table}[htb]
\vskip 1cm
\caption{Simple Jordan pairs}\label{simple-jordan}

\medskip
\begin{minipage}{\linewidth}
\begin{tabular}{lllll}
\hline

& & & &\\[-2ex]
{\scriptsize $\mathcal{U}$-Type} & {\scriptsize $(\mathcal{U}^+,\mathcal{U}^-)$-description} & {\scriptsize  $\Inder\, \mathcal{U}$} &
{\scriptsize  $\mathcal{U}^+$} & {\scriptsize  $\mathcal{U}^-$}\\

& & & &\\[-2ex]
\hline
&&&&\\[-2ex]

 {\scriptsize  $\mathrm{I}_{p,q}\footnote{The isomorphism $\mathrm{I}_{2,2}\cong \mathrm{IV}_2$ has been omitted in the classification given in \cite{Loos}.} $} & {\scriptsize  $\mathcal{U}^+={\cM}_{p,q}(k)$} & {\scriptsize  $A_{p-1}\tiny{\times} A_{q-1}\tiny{\times} Z$} &
{\scriptsize $V(\lambda_1)\tiny{\otimes} V(\lambda'_1)$} & {\scriptsize $V(\lambda_{p-1})\tiny{\otimes} V(\lambda'_{q-1})$}\\[-0,5ex]
{\tiny  $p\tiny{\geq} q \tiny{\geq} 1$}&{\scriptsize  $\mathcal{U}^-={\cM}_{p,q}(k)$} & & &\\[-0,5ex]
&{\tiny  $\{xyz\}=xy^tz+zy^tx$}& & &\\

& & & &\\[-2ex]

 {\scriptsize  $\mathrm{II}_n $} & {\scriptsize  $\mathcal{U}^+=\mathcal{U}^-=\mathrm{\cA}_n(k)$} & {\scriptsize  $A_{n-1}\times Z$} &
{\scriptsize $V(\lambda_2)$} & {\scriptsize $V(\lambda_{n-2})$}\\[-0,5ex]
{\tiny  $n\geq 5$}&{\tiny  $\{xyz\}=xy^tz+zy^tx$}& & &\\
& & & &\\[-2ex]
 {\scriptsize  $\mathrm{III}_n $} & {\scriptsize  $\mathcal{U}^+=\mathcal{U}^-=\cH_n(k)$} & {\scriptsize  $A_{n-1}\times Z$} &
{\scriptsize $V(2\lambda_1)$} & {\scriptsize $V(2\lambda_{n-1})$}\\[-0,5ex]
{\tiny  $n\geq 2$}&{\tiny  $\{xyz\}=xy^tz+zy^tx$}& & &\\
& & & &\\[-2ex]

 {\scriptsize  $\mathrm{IV}_{2n} $} & {\scriptsize  $\mathcal{U}^+=\mathcal{U}^-=k^{2n}$} & {\scriptsize  $D_{n}\times Z$} &
{\scriptsize $V(\lambda_1)$} & {\scriptsize $V(\lambda_{1})$}\\[-0,5ex]
{\tiny  $n\geq 3$}&{\tiny $\{xyz\}\tiny{=}b(x\tiny{,}y)z\tiny{+}b_{x,y}(z)$} & & &\\[-0,5ex]
&{\tiny $b(x,y)=b(y,x)$}& & &\\

& & & &\\[-2ex]

{\scriptsize  $\mathrm{IV}_{2n+1} $} & {\scriptsize  $\mathcal{U}^+=\mathcal{U}^-=k^{2n+1}$} & {\scriptsize  $B_{n}\times Z$} &
{\scriptsize $V(\lambda_1)$} & {\scriptsize $V(\lambda_{1})$}\\[-0,5ex]
{\tiny  $n\geq 2$}&{\tiny $\{xyz\}\tiny{=}b(x\tiny{,}y)z\tiny{+}b_{x,y}(z)$}& & &\\[-0,5ex]
&{\tiny $b(x,y)=b(y,x)$}& & &\\
& & & &\\[-2ex]

{\scriptsize  $\mathrm{V}$} & {\scriptsize  $\mathcal{U}^+=\mathcal{U}^-={\cM}_{1,2}(\cO)$} & {\scriptsize  $D_{5}\times Z$} &
{\scriptsize $V(\lambda_4)$} & {\scriptsize $V(\lambda_{5})$}\\[-0,5ex]
&{\tiny  $\{xyz\}=x\bar{y}^tz+z\bar{y}^tx$}& & &\\
& & & &\\[-2ex]

{\scriptsize  $\mathrm{VI}$} & {\scriptsize  $\mathcal{U}^+=\mathcal{U}^-=\cH_3(\cO)$} & {\scriptsize  $E_6\times Z$} &
{\scriptsize $V(\lambda_1)$} & {\scriptsize $V(\lambda_{6})$}\\[-0,5ex]
&{\tiny  $\{xyz\}\tiny{=}x(zy)\tiny{+}z(xy)\tiny{-}(zx)y$} & & &\\
& & & &\\[-1ex]

\hline
\end{tabular}
\end{minipage}
\end{table}

\begin{table}[htb]
\vskip 1cm
\caption{Simple anti-Jordan pairs}\label{simple-antijordan}
\medskip
\begin{minipage}{\linewidth}
\begin{tabular}{lllll}
\hline
& & & &\\[-2ex]
{\scriptsize $\mathcal{U}$-Type} & {\scriptsize $(\mathcal{U}^+,\mathcal{U}^-)$-description} & {\scriptsize  $\Inder\, \mathcal{U}$} &
{\scriptsize  $\mathcal{U}^+$} & {\scriptsize  $\mathcal{U}^-$}\\

& & & &\\[-2ex]
\hline
&&&&\\[-2ex]

 {\scriptsize  $\mathrm{GL}_{p,q} $} & {\scriptsize  $\mathcal{U}^+={\cM}_{p,q}(k)$}  & {\scriptsize  $A_{p-1}\tiny{\times} A_{q-1}\tiny{\times} Z\  {\tiny p\neq q}$} &
{\scriptsize $V(\lambda_1)\tiny{\otimes} V(\lambda'_1)$} & {\scriptsize $V(\lambda_{p-1})\tiny{\otimes} V(\lambda'_{q-1})$}\\[-0,5ex]
{\tiny  $p\tiny{\geq} q \tiny{\geq} 1$}&{\scriptsize  $\mathcal{U}^-={\cM}_{p,q}(k)$} & {\scriptsize $A_{p-1}\tiny{\times} A_{q-1}\quad \tiny{p= q}$}& &\\[-0,5ex]
{\tiny  $pq > 1$}&{\tiny  $\{xyz\}=xy^tz-zy^tx$} & & &\\
& & & &\\[-1ex]

 {\scriptsize  $\mathrm{Sps}(2n)$} & {\scriptsize  $\mathcal{U}^+=\mathcal{U}^-=k^{2n}$} & {\scriptsize  $C_{n}\times Z$} &
{\scriptsize $V(\lambda_1)$} & {\scriptsize $V(\lambda_{1})$}\\[-0,5ex]
{\tiny  $n\geq 1$}&{\tiny  $\{xyz\}\tiny{=}b(x\tiny{,}y)z\tiny{+}b_{x,y}(z)$} & & &\\[-0,5ex]
&{\tiny  $b(x,y)=-b(y,x)$} & & &\\
& & & &\\[-1ex]

 {\scriptsize   $\mathrm{Sym}(n)$}\footnote{\ In \cite[Proposition 2.8]{FauFe} the anti-Jordan pair $\mathrm{Sym}(2)$ is erroneously included as simple: $\{\cA_2(k)\cH_2(k)\cA_2(k)\}=0$, so $(\cH_2(k),0)$ is a proper ideal.} & {\scriptsize  $\mathcal{U}^+=\cH_n(k)$} & {\scriptsize  $A_{n-1}$} &
{\scriptsize $V(2\lambda_1)$} &{\scriptsize $V(\lambda_{n-2})$} \\[-0,5ex]
{\tiny  $n\geq 3$}&{\scriptsize  $\mathcal{U}^-=\mathrm{\cA}_n(k)$} & & &\\[-0,5ex]
&{\tiny  $\{xyz\}=xyz-zyx$} & & &\\
& & & &\\[-2ex]
\hline
\end{tabular}
\end{minipage}
\end{table}

\clearpage

\section{Epilogue}

The aim of this final section is to sumarize the complete classification of irreducible LY-algebras while emphasizing their connections to other algebraic systems. Following \cite[Theorem 2.4]{BEM09} we arrive at the irreducible LY-algebras of Adjoint Type. They are nothing else but simple Lie algebras with binary and ternary products given by the Lie bracket as Table~\ref{clasificacion adjuntos} shows. From \cite[Theorems 4.1 and 4.4]{BEM09} we get the information given in Table~\ref{clasificacion no-simple}. In this table, the irreducible  LY-algebras of non-simple type and exceptional enveloping algebra appear related to the Classical Tits Construction and symplectic (equivalently Freudenthal) triple systems $\mathcal{T}_\mathcal{J}$ attached to a Jordan simple algebra $\mathcal{J}$ of degree 3 or equal to the base field $k$. In the classical enveloping algebra case, the non-simple classification follows from a slight generalization  of the Tits Construction due to G. Benkart and E. Zelmanov and given in \cite{BenZel}. Along this paper we have seen that in the generic case, apart from the Lie triple systems and the exceptional cases $(G_2, \spl_2(k))$, $(E_7, G_2)$ and $(E_7, \spl_3(k))$, the irreducible LY-algebras are related to reductive pairs $(\spl(V)\  \mathrm{or}\  \so(V)\ \mathrm{or}\  \frsp(V), \Der^\star V)$ for a suitable triple system $V$ with $\Der^\star V$ closely related to the (inner) derivation algebra of the system. In this way, either Jordan or anti-Jordan pairs (triple system) appear in the $\spl$-case, Lie or orthogonal triple systems in the $\so$-case and symplectic or anti-Lie triple systems in the $\frsp$-case. This yields  our final Table \ref{clasificacion generico} according to Theorems \ref{descripcion-sl}, \ref{descripcion-so} and \ref{teo-sp}. Note that, apart from simple Lie algebras, the basic ingredients in the classification are the composition algebras ($k,k\times k, \mathcal{Q}$ and $\mathcal{O}$), and simple Jordan algebras and their zero trace elements ($\mathcal{H}_n(k)), \mathcal{H}_n(k\times k), \mathcal{H}_n(\mathcal{Q})$, $\mathcal{H}_n(\mathcal{O})$ and $\mathcal{J}(k^n)=k1\oplus k^n$ the Jordan algebra of a nondegenerate symmetric bilinear form, so $\mathcal{J}(k^n)_0=k^n$).

\begin{table}[htb]
\vskip 1cm
\caption{Irreducible LY-algebras of Adjoint Type}\label{clasificacion adjuntos}

\medskip
\begin{minipage}{\linewidth}
\begin{center}
\begin{tabular}{lllll}
\hline

 & &\\[-2ex]
{\scriptsize $\mathfrak{g(\m)}$} & &{\scriptsize $\mathfrak{h}$} & &{\scriptsize $\mathfrak{m}$-description}\\

 & &\\[-2ex]
\hline
&&\\[-2ex]
 {\scriptsize  $k[t]/(t^2-1)\otimes \mathcal{L}$} & &{\scriptsize  $\Der \mathcal{L}\cong  \mathcal{L}$}&  & {\scriptsize  $\mathcal{L}$\footnote{\ $\mathcal{L}$ stands for a simple Lie algebra with product $[a,b]$, so $\mathcal{L}$ is either a classical linear algebra $\spl_n(k), n\geq 1$ (Cartan type $A_{n-1}$), $\so_n(k), n\geq 5$ (Cartan type $B_{k}$ or $D_k$ according to $n=2k+1$ or $n=2k$), $\frsp_{2n}, n\geq 3$ (Cartan type $C_{n}$) or an exceptional algebra of type $G_2, F_4, E_6,E_7,E_8$.}}\\[-0,5ex]
& &&&{\scriptsize $a\cdot b=0$}\\[-0,5ex]
& &&&{\scriptsize  $[a,b,c]=[[a,b],c]$}\\[1ex]

{\scriptsize  $k[t]/(t^2-t-\beta)\otimes \mathcal{L}$} & &{\scriptsize  $\Der \mathcal{L}\cong  \mathcal{L}$}&
& {\scriptsize  $\mathcal{L}$}\\[-0,5ex]
{\tiny $\beta \neq -1/4$}& &&&{\scriptsize $a\cdot b=[a,b]$}\\[-0,5ex]
&&&&{\scriptsize  $[a,b,c]=\beta[[a,b],c]$}\\[1ex]

{\scriptsize  $k[t]/(t^2)\otimes \mathcal{L}$}& &{\scriptsize  $\Der \mathcal{L}\cong  \mathcal{L}$}&
& {\scriptsize  $\mathcal{L}$}\\[-0,5ex]
 &&&&{\scriptsize $a\cdot b=[a,b]$}  \\[-0,5ex]
&&&&{\scriptsize $[a,b,c]=-1/4[[a,b],c]$}\\[1ex]

&&& &\\[-1ex]

\hline
\end{tabular}
\end{center}
\end{minipage}
\end{table}

\begin{table}[htb]
\vskip 1cm
\caption{Irreducible LY-algebras of Non-simple Type}\label{clasificacion no-simple}

\medskip
\begin{minipage}{\linewidth}
\begin{tabular}{lllll}
\hline

& & & &\\[-2ex]
 {\scriptsize $\mathfrak{g(\m)}$} & &{\scriptsize $\mathfrak{h}$} &
 & {\scriptsize $\mathfrak{m}$-description}\\

& & & &\\[-2ex]
\hline
&&&&\\[-2ex]

{\scriptsize  $\spl_{pq}(k)$} & &{\scriptsize  $\spl_{p}(k) \oplus \spl_{q}(k)$} &
& {\scriptsize  $\spl_p(k)\otimes \spl_q(k)$}\\[-0,5ex]
{\scriptsize $2\leq p\leq q$}  & & & &{\scriptsize $(a\otimes f)\cdot(b\otimes g)=\frac{1}{2}[a,b]\otimes
(fg+gf-\frac{2}{q}\tr(fg)I_{q})+
$}\\[-0,5ex]

{\scriptsize $(p,q)\neq (2,2)$}  & & & &{\scriptsize $
\frac{1}{2}(ab+ba-\frac{2}{p}\tr(ab)I_{p})\otimes [f,g] $}\\[-0,5ex]
& && &{\scriptsize  $[a\otimes f,b\otimes g,c\otimes h]=
  \frac{1}{q}[[a,b],c]\otimes \tr(fg)h+$}\\[-0,5ex]

& && &{\scriptsize  $\frac{1}{p}\tr(ab)c\otimes
  [[f,g],h]$}\\[1ex]

{\scriptsize  $\so_{p+q}(k)$}&&{\scriptsize  $\so_p(k)\oplus \so_q(k)$}  & &{\scriptsize  $k^p\otimes k^q$}\\[-0,5ex]
{\scriptsize  $3\leq p\leq q$}& & & &{\scriptsize   $(u\otimes x) \cdot (v\otimes y)=0$}\\[-0,5ex]
& & & &{\scriptsize   $[u\otimes x, v\otimes y,w\otimes z]=\varphi(x,y)(\psi_{u,v}(w)\otimes
 z)+$}\\[-0,5ex]
& & & &{\scriptsize  $\psi(u,v)(w\otimes
\varphi_{x,y}(z))$}\\[-0,5ex]
&&&&{\scriptsize  $b=\varphi, \psi:\,   b(x,y)=b(y,x)$ and $b_{x,y}(z)=b(x,z)y-b(y,z)x$}\\[1ex]

 {\scriptsize  $\so_{4q}(k)$} & &{\scriptsize  $\frsp_{2}(k) \oplus \frsp_q(k)$} &
& {\scriptsize  $\frsp_2(k)\otimes \mathcal{H}_q(\mathcal{Q})_{_0}$}\\[-0,5ex]
{\scriptsize $3\leq q$}  & & & &{\scriptsize $(a\otimes f)\cdot(b\otimes g)=\frac{1}{2}[a,b]\otimes
(fg+gf-\frac{2}{q}\tr(fg)I_{q})
$}\\[-0,5ex]

& && &{\scriptsize  $[a\otimes f,b\otimes g,c\otimes h]=
  \frac{1}{q}[[a,b],c]\otimes\tr(fg)h+$}\\[-0,5ex]

& && &{\scriptsize  $\frac{1}{2}\tr(ab)c\otimes
  [[f,g],h]$}\\[1ex]

{\scriptsize  $\frsp_{p+q}(k)$}&&{\scriptsize  $\frsp_p(k)\oplus \frsp_q(k)$}  & &{\scriptsize  $k^p\otimes k^q$}\\[-0,5ex]
{\scriptsize  $2\leq p\leq q$}& & & &{\scriptsize   $(u\otimes x) \cdot (v\otimes y)=0$}\\[-0,5ex]
&&&&{\scriptsize   $[u\otimes x, v\otimes y,w\otimes z]=\varphi(x,y)(\psi_{u,v}(w)\otimes
 z)+$}\\[-0,5ex]
& & & &{\scriptsize  $\psi(u,v)(w\otimes
\varphi_{x,y}(z))$}\\[-0,5ex]
&&&&{\scriptsize  $b=\varphi, \psi:\,   b(x,y)=-b(y,x)$ and $b_{x,y}(z)=b(x,z)y+b(y,z)x$}\\[1ex]

{\scriptsize  $\frsp_{2q}(k)$} & &{\scriptsize  $\frsp_{2}(k) \oplus \so_{q}(k)$} &
& {\scriptsize  $\frsp_2(k)\otimes \mathcal{H}_q(k)_{_0}$}\\[-0,5ex]
{\scriptsize $3\leq q$}  & & & &{\scriptsize $(a\otimes f)\cdot(b\otimes g)=\frac{1}{2}[a,b]\otimes
(fg+gf-\frac{2}{q}\tr(fg)I_{q})
$}\\[-0,5ex]

& && &{\scriptsize  $[a\otimes f,b\otimes g,c\otimes h]=
 \frac{1}{q} [[a,b],c]\otimes\tr(fg)h+$}\\[-0,5ex]

& && &{\scriptsize  $\frac{1}{2}\tr(ab)c\otimes
  [[f,g],h]$}\\[1ex]

 {\scriptsize  $G_2$}  & &{\scriptsize  $\frsp_2(k)\oplus \spl_2(k)$} & &{\scriptsize  $k^2\otimes \mathcal{T}_k$}\\[-0,5ex]
 & & & &{\scriptsize  $(u\otimes x)\cdot(v\otimes y)=0$}\\[-0,5ex]
       &&&&{\scriptsize   $[u\otimes x, v\otimes y,w\otimes z]=\varphi(x,y)(\psi_{u,v}(w)\otimes
 z)+\psi(u,v)w \otimes xyz$}\\[-0,5ex]
&&&&{\scriptsize  $b=\varphi, \psi:\,   b(x,y)=-b(y,x)$ and $b_{x,y}(z)=b(x,z)y+b(y,z)x$}\\[1ex]

   {\scriptsize  $F_4$}  & &{\scriptsize  $\frsp_2(k)\oplus \frsp_6(k)$} & &{\scriptsize  $k^2\otimes \mathcal{T}_{\mathcal{H}_3(k)}$}\\[-0,5ex]
 & & & &{\scriptsize  $(u\otimes x)\cdot(v\otimes y)=0$}\\[-0,5ex]
&&&&{\scriptsize   $[u\otimes x, v\otimes y,w\otimes z]=\varphi(x,y)(\psi_{u,v}(w)\otimes
 z)+\psi(u,v)w \otimes xyz$}\\[-0,5ex]
&&&&{\scriptsize  $b=\varphi, \psi:\,   b(x,y)=-b(y,x)$ and $b_{x,y}(z)=b(x,z)y+b(y,z)x$}\\[1ex]

{\scriptsize  $F_4$}  & &{\scriptsize  $G_2\oplus \spl_2(k)$} & &{\scriptsize  $\mathcal{O}_0\otimes \mathcal{H}_3(k)_{_0}$}\\[-0,5ex]
& & & &{\scriptsize  $(a\otimes x)\cdot(b\otimes y)=\frac{1}{2}[a,b]\otimes (x\bullet y-t(x\bullet
y)1)$}\\[-0,5ex]
& & & &{\scriptsize  $[a\otimes x,b\otimes y,c\otimes z]=
  D_{a,b}(c)\otimes t(x\bullet y)z +t(ab)c\otimes
  d_{x,y}(z)$}\\[-0,5ex]
& & & &{\scriptsize $x\bullet y=\frac{1}{2}(xy+yx)$ and $d_{x,y}(z)=x\bullet(y\bullet z)-y\bullet(x\bullet z)$}\\[-0,5ex]
& & & &{\scriptsize $D_{a,b}(c)=\frac{1}{4}([[a,b],c]+3((ac)b-a(cb))$}\\[-0,5ex]
& & & &{\scriptsize $t(ab)$ and $t(x\bullet
y)$ the normalized traces}\\[1ex]

 {\scriptsize  $E_6$}  & &{\scriptsize  $\frsp_2(k)\oplus \spl_6(k)$} & &{\scriptsize  $k^2\otimes \mathcal{T}_{\mathcal{H}_3(\mathcal{K})}$}\\[-0,5ex]
 & & & &{\scriptsize  $(u\otimes x)\cdot(v\otimes y)=0$}\\[-0,5ex]
 &&&&{\scriptsize   $[u\otimes x, v\otimes y,w\otimes z]=\varphi(x,y)(\psi_{u,v}(w)\otimes
 z)+\psi(u,v)w \otimes xyz$}\\[-0,5ex]
&&&&{\scriptsize  $b=\varphi, \psi:\,   b(x,y)=-b(y,x)$ and $b_{x,y}(z)=b(x,z)y+b(y,z)x$}\\[1ex]
  {\scriptsize  $E_6$}  & &{\scriptsize  $G_2\oplus \spl_3(k)$} & &{\scriptsize  $\mathcal{O}_0\otimes \mathcal{H}_3(k\times k)_{_0}$}\\[-0,5ex]
  & & & &{\scriptsize  $(a\otimes x)\cdot(b\otimes y)=\frac{1}{2}[a,b]\otimes (x\bullet y-t(x\bullet
y)1)$}\\[-0,5ex]
& & & &{\scriptsize  $[a\otimes x,b\otimes y,c\otimes z]=
  D_{a,b}(c)\otimes t(x\bullet y)z +t(ab)c\otimes
  d_{x,y}(z)$}\\[-0,5ex]
& & & &{\scriptsize $x\bullet y=\frac{1}{2}(xy+yx)$ and $d_{x,y}(z)=x\bullet(y\bullet z)-y\bullet(x\bullet z)$}\\[-0,5ex]
& & & &{\scriptsize $D_{a,b}(c)=\frac{1}{4}([[a,b],c]+3((ac)b-a(cb))$}\\[-0,5ex]
& & & &{\scriptsize $t(ab)$ and $t(x\bullet
y)$ the normalized traces}\\[1ex]

&&& &\\[-1ex]

\end{tabular}
\end{minipage}
\end{table}

\begin{table}[htb]
\vskip 0,5cm

\medskip
\begin{minipage}{\linewidth}
\begin{tabular}{lllll}

& & & &\\[-2ex]

&&&&\\[-2ex]

{\scriptsize  $E_7$}  & &{\scriptsize  $\frsp_2(k)\oplus \so_{12}(k)$} & &{\scriptsize  $k^2\otimes \mathcal{T}_{\mathcal{H}_3(\mathcal{Q})}$}\\[-0,5ex]
 & & & &{\scriptsize  $(u\otimes x)\cdot(v\otimes y)=0$}\\[-0,5ex]
  &&&&{\scriptsize   $[u\otimes x, v\otimes y,w\otimes z]=\varphi(x,y)(\psi_{u,v}(w)\otimes
 z)+\psi(u,v)w \otimes xyz$}\\[-0,5ex]
&&&&{\scriptsize  $b=\varphi, \psi:\,   b(x,y)=-b(y,x)$ and $b_{x,y}(z)=b(x,z)y+b(y,z)x$}\\[1ex]

   {\scriptsize  $E_7$}  & &{\scriptsize  $G_2\oplus \frsp_6(k)$} & &{\scriptsize  $\mathcal{O}_0\otimes \mathcal{H}_3(\mathcal{Q})_{_0}$}\\[-0,5ex]
    & & & &{\scriptsize  $(a\otimes x)\cdot(b\otimes y)=\frac{1}{2}[a,b]\otimes (x\bullet y-t(x\bullet
y)1)$}\\[-0,5ex]
& & & &{\scriptsize  $[a\otimes x,b\otimes y,c\otimes z]=
  D_{a,b}(c)\otimes t(x\bullet y)z +t(ab)c\otimes
  d_{x,y}(z)$}\\[-0,5ex]
& & & &{\scriptsize $x\bullet y=\frac{1}{2}(xy+yx)$ and $d_{x,y}(z)=x\bullet(y\bullet z)-y\bullet(x\bullet z)$}\\[-0,5ex]
& & & &{\scriptsize $D_{a,b}(c)=\frac{1}{4}([[a,b],c]+3((ac)b-a(cb))$}\\[-0,5ex]
& & & &{\scriptsize $t(ab)$ and $t(x\bullet
y)$ the normalized traces}\\[1ex]

     {\scriptsize  $E_7$}  & &{\scriptsize  $\spl_2(k)\oplus F_4$} & &{\scriptsize  $\mathcal{Q}_0\otimes \mathcal{H}_3(\mathcal{Q})_{_0}$}\\[-0,5ex]

     & & & &{\scriptsize  $(a\otimes x)\cdot(b\otimes y)=\frac{1}{2}[a,b]\otimes (x\bullet y-t(x\bullet
y)1)$}\\[-0,5ex]
& & & &{\scriptsize  $[a\otimes x,b\otimes y,c\otimes z]=
  D_{a,b}(c)\otimes t(x\bullet y)z +t(ab)c\otimes
  d_{x,y}(z)$}\\[-0,5ex]
& & & &{\scriptsize $x\bullet y=\frac{1}{2}(xy+yx)$ and $d_{x,y}(z)=x\bullet(y\bullet z)-y\bullet(x\bullet z)$}\\[-0,5ex]
& & & &{\scriptsize $D_{a,b}(c)=\frac{1}{4}[[a,b],c]$}\\[-0,5ex]
& & & &{\scriptsize $t(ab)$ and $t(x\bullet
y)$ the normalized traces}\\[1ex]

   {\scriptsize  $E_8$}  & &{\scriptsize  $\frsp_2(k)\oplus E_7$} & &{\scriptsize  $k^2\otimes \mathcal{T}_{\mathcal{H}_3(\mathcal{O})}$}\\[-0,5ex]
  & & & &{\scriptsize  $(u\otimes x)\cdot(v\otimes y)=0$}\\[-0,5ex]
  &&&&{\scriptsize   $[u\otimes x, v\otimes y,w\otimes z]=\varphi(x,y)(\psi_{u,v}(w)\otimes
 z)+\psi(u,v)w \otimes xyz$}\\[-0,5ex]
&&&&{\scriptsize  $b=\varphi, \psi:\,   b(x,y)=-b(y,x)$ and $b_{x,y}(z)=b(x,z)y+b(y,z)x$}\\[1ex]

  {\scriptsize  $E_8$}  & &{\scriptsize  $G_2\oplus F_4$} & &{\scriptsize  $\mathcal{O}_0\otimes \mathcal{H}_3(\mathcal{O})_{_0}$}\\[-0,5ex]
   & & & &{\scriptsize  $(a\otimes x)\cdot(b\otimes y)=\frac{1}{2}[a,b]\otimes (x\bullet y-t(x\bullet
y)1)$}\\[-0,5ex]

& & & &{\scriptsize  $[a\otimes x,b\otimes y,c\otimes z]=
  D_{a,b}(c)\otimes t(x\bullet y)z +t(ab)c\otimes
  d_{x,y}(z)$}\\[-0,5ex]
& & & &{\scriptsize $x\bullet y=\frac{1}{2}(xy+yx)$ and $d_{x,y}(z)=x\bullet(y\bullet z)-y\bullet(x\bullet z)$}\\[-0,5ex]
& & & &{\scriptsize $D_{a,b}(c)=\frac{1}{4}([[a,b],c]+3((ac)b-a(cb))$}\\[-0,5ex]
& & & &{\scriptsize $t(ab)$ and $t(x\bullet
y)$ the normalized traces}\\[1ex]

&&& &\\[-1ex]

\hline
\end{tabular}
\end{minipage}
\end{table}

\begin{table}[htb]
\vskip 1cm
\caption{Irreducible LY-algebras of Generic Type}\label{clasificacion generico}

\medskip
\begin{minipage}{\linewidth}
\begin{tabular}{lllllll}
\hline

& & & & & &\\[-2ex]
 {\scriptsize $\mathfrak{g(\m)}$} &{\scriptsize $\mathfrak{h}$} &
 {\scriptsize $(\mathfrak{g(\m)},\mathfrak{h})$-pair description}& & & &{\scriptsize $\mathfrak{m}$-description}\\

& & & & & &\\[-2ex]
\hline
&&&&&&\\[-2ex]

{\scriptsize  $\spl_{n}(k)$} & {\scriptsize  $\so_{n}(k)$}& {\scriptsize  $(\mathcal{L} \mathrm{ie}_{_0}(\mathcal{H}_n(k)), \Der \mathcal{H}_n(k))^*$} & & &
& {\scriptsize  $\mathcal{H}_n(k)_{_0}$}\\[-0,5ex]
{\scriptsize $5\leq n$}  & & & & & &{\scriptsize $a\cdot b=0
$}\\[-0,5ex]
& & & & & &{\scriptsize $[a,b,c]=(bc)a-b(ac)
$}
\\[1ex]

{\scriptsize  $\spl_{2n}(k)$} & {\scriptsize  $\frsp_{2n}(k)$}&{\scriptsize  $(\mathcal{L} \mathrm{ie}_{_0}(\mathcal{H}_n(\mathcal{Q})), \Der \mathcal{H}_n(\mathcal{Q}))$}  & & &
& {\scriptsize  $\mathcal{H}_n(\mathcal{Q})_{_0}$}\\[-0,5ex]
{\scriptsize $2\leq n$}  & & & & & &{\scriptsize $a\cdot b=0
$}\\[-0,5ex]
& & & & & &{\scriptsize $[a,b,c]=(bc)a-b(ac)\qquad
$}
\\[1ex]

{\scriptsize  $\spl_{\frac{n(n+1)}{2}}(k)$} & {\scriptsize  $\spl_{n}(k)$}&{\scriptsize  $(\spl(\mathcal{H}_n(k)), \mathcal{L}_{_0}(\mathcal{H}_n(k)))^*$} & & &
& {\scriptsize  $\mathfrak{h}^\perp$}\\[-0,5ex]
{\scriptsize $2\leq n$}  & & & & & &\\[1ex]

{\scriptsize  $\spl_{\frac{n(n-1)}{2}}(k)$} & {\scriptsize  $\spl_{n}(k)$}&{\scriptsize  $(\spl(\mathcal{A}_n(k)), \mathcal{L}_{_0}(\mathcal{A}_n(k)))$} & & &
& {\scriptsize  $\mathfrak{h}^\perp$}\\[-0,5ex]
{\scriptsize $5\leq n$}  & & & & & &\\[1ex]

{\scriptsize  $\spl_{16}(k)$} & {\scriptsize  $\so_{10}(k)$}&{\scriptsize  $(\spl(\mathcal{M}_{1,2}(\mathcal{O})), \mathcal{L}_{_0}(\mathcal{M}_{1,2}(\mathcal{O})))$} & &&
& {\scriptsize  $\mathfrak{h}^\perp$}\\[1ex]

{\scriptsize  $\spl_{27}(k)$} & {\scriptsize  $E_6$}& {\scriptsize  $(\spl(\mathcal{H}_3(\mathcal{O})), \mathcal{L}_{_0}(\mathcal{H}_3(\mathcal{O})))$} & &&
&  {\scriptsize  $\mathfrak{h}^\perp$}\\[1ex]
&&&&&&\\[-1ex]

\end{tabular}
\end{minipage}
\end{table}

\clearpage

\begin{table}[htb]


\medskip
\hskip 1cm
\begin{minipage}{\linewidth}
\begin{center}
\begin{tabular}{lllll}


&&&&\\[-2ex]

{\scriptsize  $\so_{n+1}(k)$} & {\scriptsize  $\so_{n}(k)$}&{\scriptsize  $(\mathcal{L} \mathrm{ie}_{_0}(\mathcal{J}(k^n)), \Der \mathcal{J}(k^n))$\footnote{$\mathcal{L} \mathrm{ie}_{_0}(\mathcal{J})$ stands for the derived algebra of the Lie multiplication algebra attached to the Jordan algebra $\mathcal{J}$ and $\mathcal{L}_{_0}(\mathcal{T})$ is as defined in Theorem \ref{descripcion-sl} for the Jordan triple $\mathcal{T}.$}} &
& {\scriptsize  $k^{n}$}\\[-0,5ex]
{\scriptsize $5\leq n$}  & & & &{\scriptsize $x\cdot y=0
$}\\[-0,5ex]
& & & &{\scriptsize $[x,y,z]=b(x,z)y-b(y,z)x
$}\\[-0,5ex]
& & & &{\scriptsize $b(x,y)=b(y,x)
$}\\[1ex]

{\scriptsize  $\so_{\dimens \mathcal{L}}(k)$} & {\scriptsize  $\Der\mathcal{L}$}& {\scriptsize  $(\so(\mathcal{L}), \Der\mathcal{L})$\footnote{$\mathcal{L}$ stands for a simple Lie algebra different from $\spl_n(k)$.}} &
&  {\scriptsize  $\mathfrak{h}^\perp$}\\[1ex]

{\scriptsize  $\so_{\frac{n^2+n-2}{2}}(k)$} & {\scriptsize  $\so_n(k)$}& {\scriptsize  $(\so(\mathcal{H}_n(k)_{_0}), \Der\mathcal{H}_n(k))$} &
&  {\scriptsize  $\mathfrak{h}^\perp$}\\[-0,5ex]
{\scriptsize $5\leq n$}  & & & &\\[1ex]

{\scriptsize  $\so_{2n^2-n-1}(k)$} & {\scriptsize  $\frsp_{2n}(k)$}& {\scriptsize  $(\so(\mathcal{H}_n(\mathcal{Q})_{_0}), \Der\mathcal{H}_n(\mathcal{Q}))$} &
&  {\scriptsize  $\mathfrak{h}^\perp$}\\[-0,5ex]
{\scriptsize $3\leq n$}  & & & &\\[1ex]

{\scriptsize  $\so_{26}(k)$} & {\scriptsize  $F_4$}& {\scriptsize  $(\so(\mathcal{H}_3(\mathcal{O})_{_0}), \Der\mathcal{H}_3(\mathcal{O}))$} &
&  {\scriptsize  $\mathfrak{h}^\perp$}\\[1ex]

{\scriptsize  $\so_{7}(k)$} & {\scriptsize  $G_2$}&{\scriptsize  $(\so(\mathcal{O}_0), \Der \mathcal{O})$}  &
&  {\scriptsize  $\mathcal{O}_0$}\\[-0,5ex]
& & & &{\scriptsize $a\cdot b=ab-ba
$}\\[-0,5ex]
& & & &{\scriptsize $[a,b,c]=2\bigl([[a,b],c]-$}\\[-0,5ex]
& & & &{\scriptsize $\qquad \qquad \quad 3\bigl((ac)b-a(cb)\bigr)\bigr)$}\\[1ex]

{\scriptsize  $\so_{16}(k)$} & {\scriptsize  $\so_9(k)$}& {\scriptsize  $(\so(\mathcal{T}_{_{\tiny (F_4,B_4)}})
, \Inder\mathcal{T}_{_{\tiny (F_4,B_4)}})$\footnote{$\mathcal{T}_{\tiny (\mathfrak {g},\mathfrak{s})}$ stands for a simple Lie triple system attached to one of the exceptional symmetric pairs $(\mathfrak {g},\mathfrak{s})=(F_4,B_4), (E_6,C_4), (E_7,A_7)$ or $(E_8,D_8)$.}} &
&  {\scriptsize  $\mathfrak{h}^\perp$}\\[1ex]

{\scriptsize  $\so_{42}(k)$} & {\scriptsize  $\frsp_8(k)$}& {\scriptsize  $(\so(\mathcal{T}_{_{\tiny (E_6,C_4)}}), \Inder\mathcal{T}_{_{\tiny (E_6,C_4)}})$} &
&  {\scriptsize  $\mathfrak{h}^\perp$}\\[1ex]

{\scriptsize  $\so_{70}(k)$} & {\scriptsize  $\spl_8(k)$}& {\scriptsize  $(\so(\mathcal{T}_{_{\tiny (E_7,A_7)}}), \Inder\mathcal{T}_{_{\tiny (E_7,A_7)}})$} &
&  {\scriptsize  $\mathfrak{h}^\perp$}\\[1ex]

{\scriptsize  $\so_{128}(k)$} & {\scriptsize  $\so_{16}(k)$}& {\scriptsize  $(\so(\mathcal{T}_{_{\tiny (E_8,D_8)}}), \Inder\mathcal{T}_{_{\tiny (E_8,D_8)}})$} &
&  {\scriptsize  $\mathfrak{h}^\perp$}\\[1ex]

{\scriptsize  $\frsp_4(k)$} & {\scriptsize  $\spl_2(k)$}& {\scriptsize  $(\frsp(\mathcal{T}_k), \Inder\mathcal{T}_k)$\footnote{$\mathcal{T}_{\mathcal{J}}$ stands for a simple symplectic Lie triple attached to a Jordan simple algebra $\mathcal{J}=k,\mathcal{H}_3(k), \mathcal{H}_3(k\times k),\mathcal{H}_3(\mathcal{Q})$ or $\mathcal{H}_3(\mathcal{O})$.}} &
&  {\scriptsize  $\mathfrak{h}^\perp$}\\[1ex]

{\scriptsize  $\frsp_{14}(k)$} & {\scriptsize  $\frsp_6(k)$}& {\scriptsize  $(\frsp(\mathcal{T}_{\mathcal{H}_3(k)})), \Inder\mathcal{T}_{{\mathcal{H}_3(k)}})$} &
&  {\scriptsize  $\mathfrak{h}^\perp$}\\[1ex]

{\scriptsize  $\frsp_{20}(k)$} & {\scriptsize  $\spl_6(k)$}& {\scriptsize  $(\frsp(\mathcal{T}_{\mathcal{H}_3(k\times k)})), \Inder\mathcal{T}_{{\mathcal{H}_3(k\times k)}})$} &
&  {\scriptsize  $\mathfrak{h}^\perp$}\\[1ex]

{\scriptsize  $\frsp_{32}(k)$} & {\scriptsize  $\so_{12}(k)$}& {\scriptsize  $(\frsp(\mathcal{T}_{\mathcal{H}_3(\mathcal{Q})})), \Inder\mathcal{T}_{{\mathcal{H}_3(\mathcal{Q})}})$} &
&  {\scriptsize  $\mathfrak{h}^\perp$}\\[1ex]

{\scriptsize  $\frsp_{56}(k)$} & {\scriptsize  $E_7$}& {\scriptsize  $(\frsp(\mathcal{T}_{\mathcal{H}_3(\mathcal{O})})), \Inder\mathcal{T}_{{\mathcal{H}_3(\mathcal{O})}})$} &
&  {\scriptsize  $\mathfrak{h}^\perp$}\\[1ex]



{\scriptsize  $G_2$}&{\scriptsize  $\spl_2(k)$}&&& {\scriptsize  $\mathfrak{h}^\perp$}\\[1ex]

{\scriptsize  $F_4$}&{\scriptsize  $\so_9(k)$}&&&{\scriptsize $\mathcal{T}_{_{\tiny (F_4,B_4)}}$}\\[1ex]

{\scriptsize  $E_6$}&{\scriptsize  $\frsp_4(k)$}&&&{\scriptsize $\mathcal{T}_{_{\tiny (E_6,C_4)}}$}\\[1ex]
{\scriptsize  $E_6$}&{\scriptsize  $G_2$}&&& {\scriptsize  $\mathfrak{h}^\perp$}\\[1ex]

{\scriptsize  $E_6$}&{\scriptsize  $F_4$}&{\scriptsize  $(\mathcal{L} \mathrm{ie}_{_0}(\mathcal{H}_3(\mathcal{O})), \Der \mathcal{H}_3(\mathcal{O}))$} &&{\scriptsize  $\mathcal{H}_3(\mathcal{O})_0$}\\[-0,5ex]
& & & &{\scriptsize $a\cdot b=0
$}\\[-0,5ex]
& & & &{\scriptsize $[a,b,c]=(bc)a-b(ac)
$}
\\[1ex]

{\scriptsize  $E_7$}&{\scriptsize  $\spl_8(k)$}&&&{\scriptsize $\mathcal{T}_{_{\tiny (E_7,A_7)}}$}\\[1ex]

{\scriptsize  $E_7$}&{\scriptsize  $\spl_3(k)$}&&& {\scriptsize  $\mathfrak{h}^\perp$}\\[1ex]

{\scriptsize  $E_8$}&{\scriptsize  $\so_{16}(k)$}&&&{\scriptsize $\mathcal{T}_{_{\tiny (E_8,D_8)}}$}\\[1ex]
&&&&\\[-1ex]

\hline
\end{tabular}
\end{center}

\end{minipage}

\end{table}

\clearpage


\providecommand{\bysame}{\leavevmode\hbox
to3em{\hrulefill}\thinspace}
\providecommand{\MR}{\relax\ifhmode\unskip\space\fi MR }
\providecommand{\MRhref}[2]{%
  \href{http://www.ams.org/mathscinet-getitem?mr=#1}{#2}
} \providecommand{\href}[2]{#2}

\end{document}